\documentclass[12pt]{amsart}
\usepackage{amssymb}
\usepackage{amsfonts}
\usepackage{amsmath}
\usepackage{array}
\usepackage{rotating}
\usepackage{epsf}

\newtheorem{example}{Example}[section]
\newtheorem{remark}[example]{Remark}
\newtheorem{theorem}[example]{Theorem}

\newtheorem{proposition}[example]{Proposition}
\newtheorem{algorithm}[example]{Algorithm}
\newtheorem{lemma}[example]{Lemma}

\def\Proof{\noindent \it Proof -- \rm}
\def\qed{\hspace{3.5mm} \hfill \vbox{\hrule height 3pt depth 2 pt width 2mm}
\bigskip}

\def\QRP{{\rm QRP}}
\def\Hdif{{\mathcal H}^{\rm dif}}
\def\Ddif{\Delta^{\rm dif}}
\def\evt{{\rm pEv}}
\def\E{{\mathbb E}}

\def\ev{{\rm Ev}}

\def\Park{{\rm Park}}

\def\NDPF{{\rm NDPF}}

\def\park{{\bf a}}

\def\<{\langle}
\def\>{\rangle}

\def\C{\operatorname{\mathbb C}}
\def\Z{\operatorname{\mathbb Z}}

\def\SG{{\mathfrak S}}

\def\a{{\bf a}}
\def\b{{\bf b}}

\def\Sym{{\bf Sym}}

\def\dim{{\rm dim}}

\def\ch{\operatorname {ch}}
\def\chb{{\bf ch}}

\def\PF{{\rm PF}}

\def\NDPF{{\rm NDPF}}
\def\HS{{\rm H}{\mathfrak S}}

\def\shuff#1#2{\mathbin{
\hbox{\vbox{ \hbox{\vrule \hskip#2 \vrule height#1 width 0pt
}%
\hrule}%
\vbox{ \hbox{\vrule \hskip#2 \vrule height#1 width 0pt
\vrule }%
\hrule}%
}}}

\def\shuf{{\mathchoice{\shuff{7pt}{3.5pt}}%
{\shuff{6pt}{3pt}}%
{\shuff{4pt}{2pt}}%
{\shuff{3pt}{1.5pt}}}}%
\def\shuffle{\,\shuf\,}

\def\Tabvrule{\vrule width-0.4pt}       
\def\Tabhrule{\hrule \hrule height-0.4pt} 
\def\Tabstrut{\vrule height2.2ex 
                     depth0.8ex  
                     width0ex    
\relax}

\def\PasCase#1{\omit%
            $\vcenter{\hbox {\vbox to 0.4pt{}}
               \hbox{\makebox[3ex]{\Tabstrut$#1$}}}%
               \Tabvrule$}
\def\PasCasePoint{\PasCase{\cdot}}
\def\DessinCarre#1{%
    \vcenter{\hbox{}\hrule
             \hbox{\vrule\makebox[3ex]{\Tabstrut$#1$}\vrule}\Tabhrule}%
             \Tabvrule}
\def\GenRuban#1{\vcenter{\halign{&$\DessinCarre{##}$\cr#1}}\egroup}

\def\sTabvrule{\vrule width-0.4pt}
\def\sTabhrule{\hrule \hrule height-0.4pt}
\def\sTabstrut{\vrule height1.6ex depth0.6ex width0ex \relax}
\def\sDessinCarre#1{%
    \vcenter{\hbox{}\hrule
             \hbox{\vrule\makebox[2.3ex]%
                  {\sTabstrut$\scriptstyle#1$}\vrule}\sTabhrule}%
             \sTabvrule}
\def\sGenRuban#1{\vcenter{\halign{&$\sDessinCarre{##}$\cr#1}}\egroup}

\def\ruban{%
  \bgroup
  \let\ =\omit
  \let\\=\cr
  \let\x=\times
  \let\.=\PasCasePoint
  \offinterlineskip
  \GenRuban}

\def\sruban{%
  \bgroup
  \let\ =\omit
  \let\x=\times
  \let\\=\cr
  \offinterlineskip
  \sGenRuban}

\title{Noncommutative Symmetric Functions and Lagrange Inversion}

\author[J.-C.~Novelli and J.-Y.~Thibon]%
{Jean-Christophe Novelli and Jean-Yves Thibon}

\address[] {Institut Gaspard Monge, Universit\'e de Marne-la-Vall\'ee \\
5 Boulevard Descartes \\Champs-sur-Marne \\77454 Marne-la-Vall\'ee cedex 2 \\
FRANCE}
\email[Jean-Christophe Novelli]{novelli@univ-mlv.fr}
\email[Jean-Yves Thibon]{jyt@univ-mlv.fr} 
\date{}

\begin{document}

\begin{abstract}
We compute the noncommutative Frobenius characteristic of the natural
action of the $0$-Hecke algebra on parking functions, and obtain as
corollaries various forms of the noncommutative Lagrange inversion formula.
\end{abstract}

\maketitle

\section{Introduction}

There are some advantages to interpret the classical Lagrange inversion
formula for the reversion of formal power series in terms of symmetric
functions (see, e.g., \cite{Mcd}, Ex. 24 p. 35, Ex. 25 p. 132, \cite{Las}
Section~2.4 and \cite{Len}). 
Recall that one possible formulation of the problem is
as follows. Given 
\begin{equation}
\varphi(x)=\sum_{n\ge 0}\varphi_nx^n\quad (\varphi_0\not=0)\,
\end{equation}
find the coefficients $c_n$ of the unique power series 
\begin{equation}
u(t)=\sum_{n\geq0}c_nt^{n+1}
\end{equation} 
satifying
\begin{equation}
t = \frac{u}{\varphi(u)} \,.
\end{equation}
We can assume without loss of generality that $\varphi_0=1$ and that
\begin{equation}
\varphi(u)=\sum_{n\ge 0}h_n(X)u^n= \prod_{n\ge 1}(1-ux_n)^{-1}=: \sigma_u(X)
\end{equation}
is the generating series of the homogeneous symmetric functions of an
infinite set of variables $X$. Indeed, the $h_n(X)$ are algebraically
independent, so that $\sigma_u(X)$ is a generic power series.

Now, symmetric functions encode various mathematical objects,
and the solution can be interpreted
in many ways, for example in terms of characters of the symmetric group.
Indeed, in the $\lambda$-ring notation, the solution reads
\begin{equation}
c_n =\frac1{n+1}h_n((n+1)X)
\end{equation}
(recall that $\sigma_t(nX)=\sigma_t(X)^n$, see, e.g., \cite{Mcd} p. 25).
On this expression, it is obvious that $c_n$ is Schur positive, in fact,
even a positive sum of homogeneous products $h_\mu$, so that it is
the Frobenius characteristic of a permutation representation of $\SG_n$.
This representation is well-known \cite{Hai1}: it is based on the set $\PF_n$
of {\em parking functions} of length $n$ (see below for the definition).
The first terms are
\begin{equation}
\label{c01234}
\begin{split}
& c_0 =1,\qquad c_1 = h_1,\qquad c_2 = h_2 + h_{11}\,, \\
& c_3 = h_{3} + 3h_{21} + h_{111}\,, \\
& c_4 = h_4 + 4h_{31} + 2h_{22} + 6h_{211} + h_{1111}\,.
\end{split}
\end{equation}
Now, we have at our disposal noncommutative analogs of the Lagrange inversion
formula~\cite{Ges,PPR,BFK}, and a theory of noncommutative symmetric functions
\cite{NCSF1,NCSF2}, known to be related to $0$-Hecke algebras in the same
way as ordinary symmetric functions are related to symmetric groups
\cite{NCSF4}. The aim of this note is to clarify the relations between these
different topics. 
We shall first analyze the natural representation of the $0$-Hecke algebra
on parking functions. This is a projective module, whose $q$-characteristic
noncommutative symmetric function turns out to be the term of degree $n$ in
the noncommutative $q$-Lagrange inversion formula.
This allows us to give simple and unified proof of all versions of the
noncommutative Lagrange formula \cite{Ges,PPR,BFK}. Interpreting the terms as
ordered trees leads to closed expressions for the expansion of the solution in
various bases.  These calculations suggest the introduction of noncommutative
analogs of Abel's polynomials, and of an infinite family of combinatorial
triangles, which includes classical refinements of the Motzkin, Catalan and
Schr\"oder numbers as the first three cases. The action of the $0$-Hecke
algebra on $(k,l)$-parking functions is also described.

\medskip
{\footnotesize
{\it Acknowledgements.-}
This project has been partially supported by CNRS and by EC's IHRP Programme,
grant HPRN-CT-2001-00272, ``Algebraic Combinatorics in Europe".
The authors would also like to thank the contributors of the MuPAD project,
and especially of the combinat part, for providing the development environment
for their research (see~\cite{HTm} for an introduction to MuPAD-Combinat).
}

\section{Notations}

Our notations for noncommutative symmetric functions will be as in
\cite{NCSF1,NCSF2}. We recall that the Hopf algebra of noncommutative symmetric
functions is denoted by $\Sym$, or by $\Sym(A)$ if we consider the realization
in terms of an auxiliary alphabet. Bases of $\Sym_n$ are labelled by
compositions $I$ of $n$. The noncommutative complete and elementary functions
are denoted by $S_n$ and $\Lambda_n$, and the notation $S^I$ means
$S_{i_1}\cdots S_{i_r}$. The ribbon basis is denoted by $R_I$.
The notation $I\vDash n$ means that $I$ is a composition of $n$.
The conjugate composition is denoted by $I^\sim$.

The graded dual of $\Sym$ is $QSym$ (quasi-symmetric functions).
The dual basis of $(S^I)$ is $(M_I)$ (monomial), and that of $(R_I)$
is $(F_I)$.

The {\em evaluation} $\ev(w)$ of a word $w$ over a totally ordered alphabet
$A$ is the sequence $(|w|_a)_{a\in A}$ where $|w|_a$ is the number of
occurences of $a$ in $w$. The {\em packed evaluation} $I=\evt(w)$ is the
composition obtained by removing the zeros in $\ev(w)$.

The {\em Hecke algebra} $H_n(q)$ ($q\in\C$) is the $\C$-algebra
generated by $n-1$ elements $T_1,\ldots,T_{n-1}$ satisfying
the braid relations and $(T_i-1)(T_i+q)=0$. We are interested in the case
$q=0$, whose representation theory can be described in terms of
quasi-symmetric functions and noncommutative symmetric functions
\cite{NCSF4,NCSF6}.

The Hopf structures on $\Sym$ and $QSym$ allows one to mimic, up
to a certain extent, the $\lambda$-ring notation which is so useful
for dealing with ordinary symmetric functions (see \cite{Las}
for the commutative version and
\cite{NCSF2} for the noncommutative extension).
If $A$ is a totally ordered alphabet, the noncommutative symmetric functions
of $nA$ $(n\in\Z)$ and $[n]_qA$ (where $[n]_q=\{1<q<\cdots<q^{n-1}\}$) are
defined by
\begin{equation}
\sigma_t(nA)=\sum_{m\ge 0}t^m S_m(nA) := \sigma_t(A)^n
\end{equation}
and
\begin{equation}
\sigma_t([n]_qA) :=\sigma_t(A)\sigma_t(qA)\cdots \sigma_t(q^{n-1}A)\,.
\end{equation}
More generally, noncommutative symmetric functions can be evaluated
on any element $x$ of a $\lambda$-ring, $S_n(x)=S^n(x)$ being the $n$-th
symmetric power. Recall that $x$ is said {\em of rank one} (resp.
{\em binomial}) if $\sigma_t(x)=(1-tx)^{-1}$ (resp.
$\sigma_t(x)=(1-t)^{-x}$). The scalar $x=1$ is the only element
having both properties. We usually consider that our auxiliary variable
$t$ is of rank one, so that $\sigma_t(A)=\sigma_1(tA)$.

For each of the noncommutative formulas obtained from representations
of the $0$-Hecke algebras, we shall give the commutative specializations
to the alphabet $A=1$ ($S_n(1)=1$ for all $n$) and to the virtual
alphabet $A=\E$, defined by $\sigma_t(\E)=e^t$. This will produce
a number of (generally known) combinatorial identities, which
can now be traced back to a common source.

\section{Permutational $0$-Hecke modules}

\subsection{}
Let $[N]=\{1,\ldots,N\}$ regarded as an ordered alphabet.
There is a  right action of $H_n(q)$ on $\C [N]^n$
corresponding to the standard right action of $\SG_n$ (see \cite{NCSF4}).
If $w=a_1a_2\cdots a_n$, one sets 
$w\cdot\sigma_i = a_1\cdots a_{i+1}a_i\cdots a_n$, and
\begin{equation}
w\cdot T_i = \left\{
\begin{array}{cccl}
w  \cdot  T_i & = &  w\cdot {\sigma_i} & \qquad
\hbox{if} \ a_i < a_{i+1} \,,  \\
w  \cdot  T_i & = & q\,  w & \qquad
\hbox{if} \ a_i = a_{i+1} \,,  \\
w  \cdot  T_i & = & q\, w\cdot{\sigma_i}
+ (q-1)\, w & \qquad \hbox{if} \ a_i > a_{i+1} \ .  \\
\end{array}
\right.
\end{equation}
For $q=0$, this simplifies as
\begin{equation}
w\cdot T_i = \left\{
\begin{array}{cccl}
w  \cdot  T_i & = &  w\cdot {\sigma_i} & \qquad
\hbox{if} \ a_i < a_{i+1} \,,  \\
w  \cdot  T_i & = & 0 & \qquad
\hbox{if} \ a_i = a_{i+1} \,,  \\
w  \cdot  T_i & = & 
-w & \qquad \hbox{if} \ a_i > a_{i+1} \ .  \\
\end{array}
\right.
\end{equation}
Thus, the image of a word $w$ by an element of $H_n(0)$ is either
(up to a sign) a rearrangement of $w$ or $0$. In particular, starting
from a nondecreasing word $v$, one obtains all rearrangements of $v$.
These form the basis of a projective $H_n(0)$-module $M$ whose noncommutative
characteristic is $\chb(M)=S^I\in\Sym=\Sym(A)$ 
where $I$ is the packed evaluation of $v$ \cite{NCSF4,NCSF6}.

The characteristic of the permutation representation
$W_n(N)=\C [N]^n$ is easily seen to be
\begin{equation}
\chb(W_n(N))= \sum_{I\vDash n}M_I(N)S^I(A) = S_n(NA)
\end{equation}
by the noncommutative Cauchy identity, since the specialization
$M_I(N):=M_I(1^N)$ (1 repeated $N$ times) of the monomial quasi-symmetric
function $M_I$ is equal to the number of words of $[N]^n$ with packed
evaluation $I$.

One can do better, and keep track of the sum of the letters, a statistic
obviously preserved by the action of $\SG_n$ or $H_n(0)$. We shall
normalize it as 
\begin{equation}
\|w\| = \sum_{i=1}^n (a_i-1)\,.
\end{equation}
Then,
\begin{equation}
\sum_{w\in[N]^n,\, \evt(w)=I}q^{\|w\|}=M_I(1,q,\ldots,q^{N-1})
=M_I([N]_q)
\end{equation}
and we can write down a $q$-characteristic
\begin{equation}
\chb_q(W_n(N))= \sum_{I\vDash n}M_I([N]_q)S^I = S_n([N]_qA)\,.
\end{equation}

\subsection{}
A \emph{parking function} on $[n]=\{1,2,\ldots,n\}$ is a word
$\park=a_1a_2\cdots a_n$ of length $n$ on $[n]$ whose non-decreasing
rearrangement $\park^\uparrow=a'_1a'_2\cdots a'_n$ satisfies $a'_i\le i$ for
all $i$.
Let $\PF_n$ be the set of such words.
We are interested in the computation of $G_n(q;A):= \chb_q(\PF_n)$.
The first values are:

\begin{equation}
\label{G01234}
\begin{split}
& G_0 =1,\qquad G_1 = S_1,\qquad G_2 = S_2 + qS^{11}\,, \\
& G_3 = S^{3} + (q+q^2)S^{21} + q^2S^{12} + q^3S^{111}\,, \\
& G_4 = S^4 + (q+q^2+q^3)S^{31} + (q^2+q^4)S^{22} + q^3S^{13} +
(q^3+q^4+q^5)S^{211}\\
&\ \ \ \ \ \ \ + (q^4+q^5)S^{121} + q^5S^{112} + q^6 S^{1111}\,.
\end{split}
\end{equation}
One can  decompose the set of words $w\in [n+r]^n$ according to the length
of their maximal parking subword $p(w)$ (which may be empty, and is clearly
unique).
If $p(w)$ is of length $k$, the complementary subword can only
involve letters greater than $k+1$, and can in fact be any word
of $[k+2,n+r]^{n-k}$. Hence \cite{KW},
\begin{equation}\label{lem1}
[n+r]^n=\bigsqcup_{k=0}^n \PF_k \shuffle [k+2,n+r]^{n-k}\,.
\end{equation}

Taking the $q$-characteristic of the underlying permutational
$0$-Hecke modules, and remembering that shuffling over disjoint
alphabets amounts to inducing representations, we obtain
\begin{equation}
S_n([n+r]_qA)= \sum_{k=0}^n \chb_q(\PF_k) q^{(k+1)(n-k)} S_{n-k}([n+r-k-1]_qA)
\end{equation}
which allows us to extract the generating series
 of $G_n(q;A):= \chb_q(\PF_n)$.
Indeed, writing
\begin{equation}
(k+1)(n-k)=\scriptstyle \binom{n+1}{2}-\binom{k+1}{2}-\binom{n-k}{2}
\end{equation}
and 
\begin{equation}
F^{(r)}(x,q;A)=\sum_{n\ge 0} x^n q^{-\binom{n}{2}}S_n([n+r]_qA)
\end{equation}
we arrive at
\begin{theorem}The generating series
 of $G_n(q;A):= \chb_q(\PF_n)$ is
\begin{equation}\label{Gr}
G(x,q;A):= \sum_{n\ge 0}x^n q^{-\binom{n+1}{2}}G_n(q;A)
= F^{(r)}(xq^{-1},q;A)\ F^{(r-1)}(x,q;A)^{-1}\,.
\end{equation}
\end{theorem}
In particular, this expression is independent of $r$, a fact which
is not easily derived by mere algebraic manipulations.

We can let $r$ tend to infinity, and obtain the simpler form
\begin{equation}
G(x,q;A)=F(xq^{-1},q;A)F(x,q;A)^{-1}
\end{equation}
where 
\begin{equation}
F(x,q;A)=\sum_{n\ge 0} x^n q^{-\binom{n}{2}}S_n\left(\frac{A}{1-q}\right)\,.
\end{equation}

\begin{example}{\rm
Let us take the specialization $A=\E$, where the ``exponential alphabet''
$\E$ is defined by $\sigma_t(\E)=e^t$ (that is, $S_n(\E)=1/n!$).
Then
\begin{equation}
\sigma_t\left(\frac{\E}{1-q}\right)=\exp\left(\frac{t}{1-q}\right)
\end{equation}
and we recover Gessel's formula for the sum enumerator of parking
functions (\cite{Ges}, see also
\cite{Stan2}, Ex. 5.48.b and 5.49.c pp. 94-95).
}
\end{example}

\begin{example}{\rm
If we take $A=1$, so that $\sigma_t(1)=(1-t)^{-1}$,
we have 
\begin{equation}
S_n\left(\frac1{1-q}\right)=\frac1{(q)_n}
\end{equation}
and replacing $q$ by $1/q$ and $x$ by $-1$
we recognize in $F(-1,1/q;1)$ and $F(-q,1/q;1)$ the left-hand sides
of the Rogers-Ramanujan identities.
We have in fact an infinity of different
expressions of the Ramanujan function $F(-qx,1/q;1)F(-x,1/q;1)^{-1}$
as $F^{(r)}(-qx,1/q;1)F^{(r-1)}(-x,1/q;1)^{-1}$. The case $r=1$ is
obtained in \cite{PPR} (precisely as an application of  noncommutative
Lagrange inversion).
}
\end{example}

\section{The functional equation}

We shall now see that $G(x,q;A)$ solves a functional equation, and
recover the noncommutative $q$-Lagrange formula in this way.
For later convenience, let us first change $q$ into $1/q$ and
consider
\begin{equation}
H(x,q;A):=G(x,q^{-1};A)=E(qx)E(x)^{-1} 
\end{equation}
where
\begin{equation}
E(x)=E(x,q;A)=\sum_{n\ge 0}x^n q^{\binom{n}{2}}
S_n\left(\frac{A}{1-q^{-1}}\right)\,.
\end{equation}
Then, $H(q^{k-1}x)=E(q^kx)E(q^{k-1}x)^{-1}$, so that
\begin{equation}
H^{(n)}(x):=H(q^{n-1}x)H(q^{n-2}x)\cdots H(x) =E(q^nx)E(x)^{-1}\,,
\end{equation}
and we can write
\begin{equation*}
\begin{split}
\sum_{n\ge 0}x^nq^{\binom{n+1}{2}}S_n(A)H^{(n)}(x)
&= \sum_{n\ge 0}x^nq^{\binom{n+1}{2}}S_n(A)E(q^nx)E(x)^{-1}\\
&=\sum_{n\ge 0}x^nq^{\binom{n+1}{2}}S_n(A)\sum_{m\ge 0}q^{\binom{m}{2}}
S_m\left(\frac{A}{1-q^{-1}}\right)E(x)^{-1}\\
&=\sum_{N\ge 0}x^Nq^{\binom{N}{2}}S_N\left(qA+\frac{A}{1-q^{-1}}\right)E(x)^{-1}\\
&=\sum_{N\ge 0}x^Nq^{\binom{N+1}{2}}S_N\left(\frac{A}{1-q^{-1}}\right)E(x)^{-1}\\
&=E(qx)E(x)^{-1}=H(x)\,.
\end{split}
\end{equation*}
The powers of $q$ can be absorbed in the products if we set $K(x)=xqH(x)$.
Finally, we obtain

\begin{theorem}
The series $K(x)=K(x,q;A)=xqG(x,q^{-1};A)$ solves the functional equation
of the noncommutative $q$-Lagrange formula of  \cite{Ges,PPR}
\begin{equation}
K(x)= qx\sum_{n\ge 0}S_n(A)\cdot K^{(n)}(x)\,.
\end{equation}
\end{theorem}
One has
\begin{equation}
\label{K01234}
\begin{split}
K(x)=& xq + x^2 q^2S_1 + x^3 (q^4S_2+q^3S^{11}) +
      x^4 (q^7S_3 + (q^5+q^6)S^{21} + q^5S^{12} + q^4S^{111})\\
& +
      x^5 (q^{11}S_4 + (q^8+q^9+q^{10})S^{31} + (q^7+q^9)S^{22} +
           q^8S^{13} + (q^6+q^7+q^8) S^{211}\\
&\ \ \ \ + (q^6+q^7)S^{121} + q^6S^{112} + q^5S^{1111})+\cdots
\end{split}
\end{equation} 
In particular, $g(A)=G(1,1;A)=\sum\chb(\PF_n)$ is the unique solution of
\begin{equation}
\label{expr-g}
g = \sum_{n\ge 0}S_n g^n\,,
\end{equation}
with $S_0=1$. The first terms are
\begin{equation}
\label{g01234}
\begin{split}
& g_0 =1,\qquad g_1 = S_1,\qquad g_2 = S_2 + S^{11}\,, \\
& g_3 = S_3 + 2S^{21} + S^{12} + S^{111}\,, \\
& g_4 = S_4 + 3S^{31} + 2S^{22} + 3S^{13} + 3S^{211} + 2S^{121} + S^{112} +
S^{1111}\,.
\end{split}
\end{equation}
Note that $g_i$ is obtained by setting
 $q=1$ in (\ref{G01234}), that (\ref{expr-g}) is 
(\ref{K01234}) for $q=x=1$ and that one recovers (\ref{c01234}) by
assuming that the $S_i$ commute.

The solution of \cite{Ges,PPR} is obtained by taking $r=1$ in
Formula~(\ref{Gr}).
The commutative image gives various forms of the Garsia-Gessel
$q$-Lagrange formula. 

\section{The general noncommutative Lagrange inversion formula}

\subsection{Nondecreasing parking functions}

The versions of \cite{Ges} and \cite{PPR} on the noncommutative
inversion formula deal with the slightly more general functional
equation
\begin{equation}
\label{genfeq}
f = S_0 +S_1 f + S_2f^2+S_3f^3+\cdots\,,
\end{equation}
where $S_0$ is another indeterminate which does not necessarily commute
with the other ones. The solution can be expressed in the form
\begin{equation}\label{gensol}
f_n = \sum_{\pi\in\NDPF_n} S^{\ev(\pi)\cdot 0}\,,
\end{equation}
where $\NDPF_n$ denotes the set of nondecreasing parking functions
on $[n]$.
For example,
\begin{eqnarray}
\label{sol1110}
f_0=S_0\,,\ f_1= S_1S_0=S^{10}\,,\ f_2=S^{110}+S^{200}\,\\
f_3=S^{1110}+S^{1200}+S^{2010}+S^{2100}+S^{3000}\,,
\end{eqnarray}
the nondecreasing parking functions giving $f_3$ being
(in this order) $123$, $122$, $113$, $112$, $111$.

\subsection{Dyck words}

Here is an amusing way to prove Formula~(\ref{gensol}), inspired by
one of the examples of \cite{PPR}. If we denote by $D$ the sum of all 
Dyck words (1 being the empty word)
\begin{equation}\label{eqdyck}
D(a,b)=1+ab +aabb +abab + aaabbb+ aabbab+ aababb+abaabb+ababab+\cdots
\end{equation}
which can be defined by the functional equation
\begin{equation}
D=1+aDbD\,,
\end{equation}
then, the series $f=Db$ satifies (\ref{genfeq}), with
\begin{equation}\label{Snab}
S_n=a^nb\,.
\end{equation} 
Indeed, iterating (\ref{eqdyck}), we have
\begin{equation}
\begin{split}
Db=&\ b+a(Db)(Db)=b+a[b+a(Db)(Db)](Db)\\
=&\ b+ab(Db)+aab(Db)(Db)+aaab(Db)(Db)(Db)+\cdots\,,
\end{split}
\end{equation}
so that we know the solution of (\ref{genfeq}) in this particular
case. But the particular case is generic:
$S=\{a^nb|n\ge 0\}$ is a prefix code, and the $S_n$ defined
by (\ref{Snab}) are algebraically independent. The general solution
(\ref{gensol}) is then obtained by decomposing the words of $Db$
on the code $S$.

This expression being granted, the other version of the solution
(as a quotient of series) is obtained directly from (\ref{lem1})
as above.

In \cite{PPR}, the specialization $S_n=\frac{1}{n!}a^nb$
is also considered, leading to what the authors have called noncommutative
inversion polynomials. 

\subsection{Trees}

\begin{figure}[ht]
\begin{center}
\leavevmode
\epsfxsize = 12cm
\epsffile{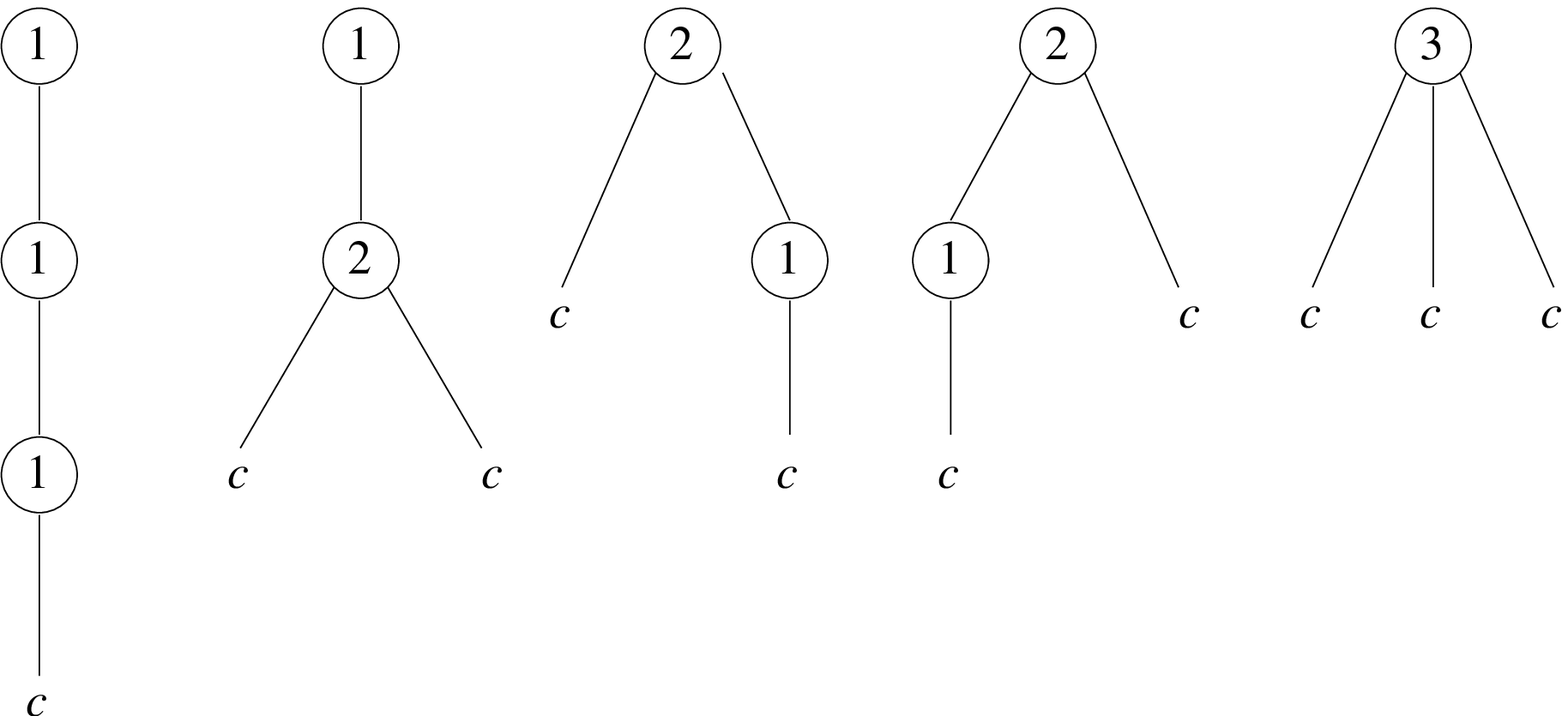}
\caption{\label{fig-5ordonnes}The five ordered trees corresponding to $f_3$.}
\end{center}
\end{figure}

Alternatively, Formula~(\ref{genfeq}) can be interpreted as a sum over ordered
trees. Let us set $c=S_0$, $d_n=S_n$, and interpret $d_n$ as the symbol of an
$n$-ary operation in Polish notation, so that for example
\begin{equation}
f_3 = d_1d_1d_1c + d_1d_2cc + d_2cd_1c + d_2d_1cc + d_3ccc
\end{equation}
is the Polish notation for
\begin{equation}
d_1(d_1(d_1(c))) + d_1(d_2(c,c)) + d_2(c,d_1(c)) + d_2(d_1(c),c) +
d_3(c,c,c)
\end{equation}
and corresponds to the five ordered trees of Figure~\ref{fig-5ordonnes}.

This implies an expression of the coefficients $\delta_I$ defined by
\begin{equation}
g_n =\ch(\PF_n)= \sum_{I\models n} \delta_I S^I
\end{equation}
since $g_n$ is obtained from $f_n$ by setting $c=1$
in (\ref{expr-g}).
Indeed,
given a tree $T$, define its \emph{skeleton} as the tree obtained by removing
the leaves $c$ and labeling the internal vertices with their arity.
Given the skeleton $S$ of a tree $T$, define its \emph{$0$-composition}
$I_0(S)$ as the sequence formed by the values of the labels of the vertices of
$S$ read in prefix order.

\begin{figure}[ht]
\begin{center}
\leavevmode
\epsfxsize =10cm
\epsffile{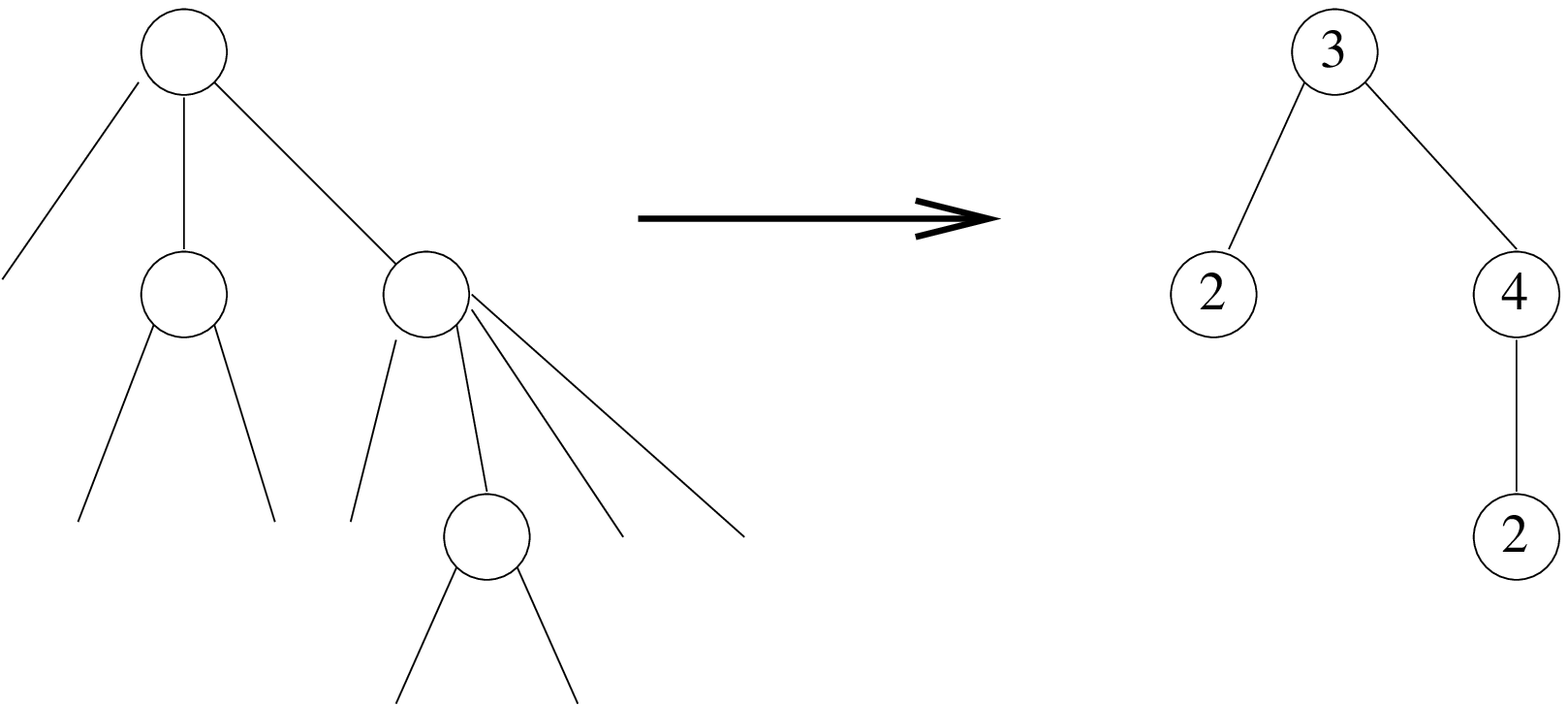}
\caption{\label{fig-squelette}A tree and its skeleton.}
\end{center}
\end{figure}

For example, one finds on Figure~\ref{fig-squelette} a tree and its skeleton.
The corresponding $0$-compo\-sition is $(3,2,4,2)$.

The number of trees with skeleton $S$ is obviously
\begin{equation}
\label{arb-squelette}
\prod_{k=1}^p \binom{i_k}{a_k}
\end{equation}
where $I_0(S)=(i_1,\ldots,i_p)$ and $a_k$ is the arity of the $k$-th vertex
of the tree $S$, numbered in prefix order.

\begin{figure}[ht]
\begin{center}
\leavevmode
\epsfxsize = 12cm
\epsffile{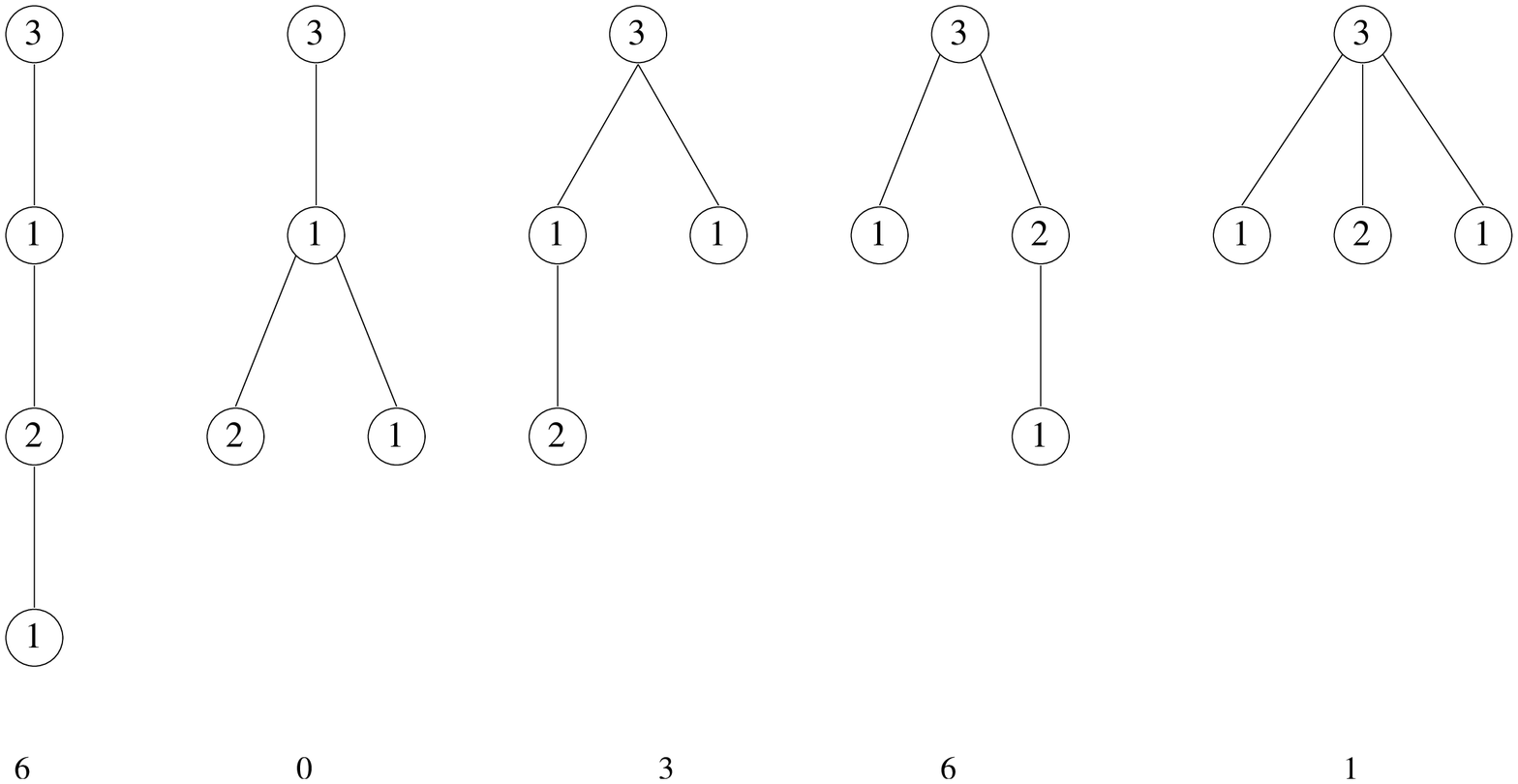}
\caption{\label{fig-ex16}Skeletons having $(3,1,2,1)$ as 
$0$-composition, and the number of trees with those skeletons.}
\end{center}
\end{figure}

For example, there are $16$ trees whose skeleton have $(3,1,2,1)$
as  $0$-composition as one can check on Figure~\ref{fig-ex16}.

Let $I=(i_1,\ldots,i_p)$ be a composition of $n$.
We are now in a position to compute $\delta_I$. Indeed, the coefficient of
$S^I$ in $g_n$ is equal to the number of ordered trees on $n+1$ vertices
whose sequence of non-zero arities in prefix order is $I$.
The skeletons of these trees are the ordered trees on $p$ vertices labeled
by the elements of $I$ in prefix order.
The sequences of arities of the skeletons in prefix order are all sequences
$(a_1,\ldots,a_p)$ such that $a_1+\cdots+a_j\geq j$ for $j<p$ and equal to
$p-1$ if $j=p$. Thus $a_p=0$ and $a_1+\cdots+a_{p-1}=p-1$ so that
\begin{equation}
\label{delta-res}
\delta_I = \sum_{(a_1,\ldots,a_{p-1})}
\prod_{k=1}^{p-1} \binom{i_k}{a_k}\,,
\end{equation}
where the sum is taken over the set of sequences $(a_1,\ldots,a_{p-1})$ such
that $a_1+\cdots+a_j\geq j$ for all $j$ and $a_1+\cdots+a_{p-1} = p-1$.

\section{Noncommutative formal diffeomorphisms}

\subsection{}
Another form of the noncommutative Lagrange inversion has been
obtained by Brouder-Frabetti-Krattenthaler \cite{BFK}. It is stated in the
form of an explicit formula for the antipode of the Hopf algebra
$\Hdif$ of ``formal diffeomorphisms''. As an associative algebra,
$\Hdif$ can be identified with $\Sym$ 
by means of the correspondence $a_n=S_n=S_n(A)$.
The coproduct can then be expressed
as
\begin{equation}
\Ddif S_n(A)=\sum_{k=0}^n S_k(A)\otimes S_{n-k}((k+1)A)\,.
\end{equation}
In  this notation, computing the antipode amounts to find a series 
\begin{equation}
h(A)=\sum_{n\ge 0} b_n 
\end{equation}
where $b_n\in\Sym_n(A)$, such that
\begin{equation}
1 = \sum_{n\ge 0} S_n(A) h(A)^{n+1}\,.
\end{equation}
Hence, $h(A)$ satisfies the functional equation
\begin{equation}\label{eqbfk}
h(A)^{-1}=\sum_{n\ge 0}S_n(A) h(A)^n\,,
\end{equation}
differing from that of Gessel and Pak-Postnikov-Retakh, which reads
\begin{equation}\label{eqppr}
g(A)=\sum_{n\ge 0}S_n(A) g(A)^n\,.
\end{equation}
However, the difference is not that big, since we have
\begin{theorem}
The relation between the noncommutative symmetric series $h(A)$ and $g(A)$
respectively defined by {\rm (\ref{eqbfk})} and {\rm (\ref{eqppr})} is
\begin{equation}
h(A)=g(-A)\,.
\end{equation}
\end{theorem}
\Proof
This is a good illustration of the power of the ``noncommutative $\lambda$-ring
notation''.
Using the expression of $g(A)$ given by putting $q=1$ in (\ref{Gr}),
we can write
\begin{equation}
g(-A)^n = F(-n)F(0)^{-1} \quad (n\in\Z)
\end{equation}
where
\begin{equation}
F(x)=\sum_{m\ge 0}S_m((x-m)A) \,.
\end{equation}
Hence,
\begin{equation*}
\begin{split}
\sum_{n\ge 0} S_n(A)g(-A)^n &= \sum_{n\ge 0} S_n(A)F(-n)F(0)^{-1}\\
&= \sum_{n\ge 0} S_n(A)\sum_{m\ge 0}S_m((-m-n)A)F(0)^{-1}\\
&= \sum_{N\ge 0}\sum_{m+n=N}S_n(A)S_m((-m-n)A)F(0)^{-1}\\
&= \sum_{N\ge 0}S_N((-N+1)A) = F(1)F(0)^{-1}=g(-A)^{-1}\,. 
\end{split}
\end{equation*}
\qed

\begin{remark}{\rm This calculation works as well for the $q$-analog,
and allows one to compute the antipode of the $q$-deformed 
coproduct
\begin{equation}
\Ddif_q S_n(A) = \sum_{k=0}^n S_n(A)\otimes S_{n-k}([k]_qA)\,.
\end{equation}
For $q=0$, this is the usual coproduct of $\Sym$. We have therefore
an interpolation between the two structures, the combinatorics
being governed by the $q$-Lagrange formula, hence by parking functions.
}
\end{remark}

\subsection{Trees}

\begin{figure}[ht]
\begin{center}
\leavevmode
\epsfxsize =14cm
\epsffile{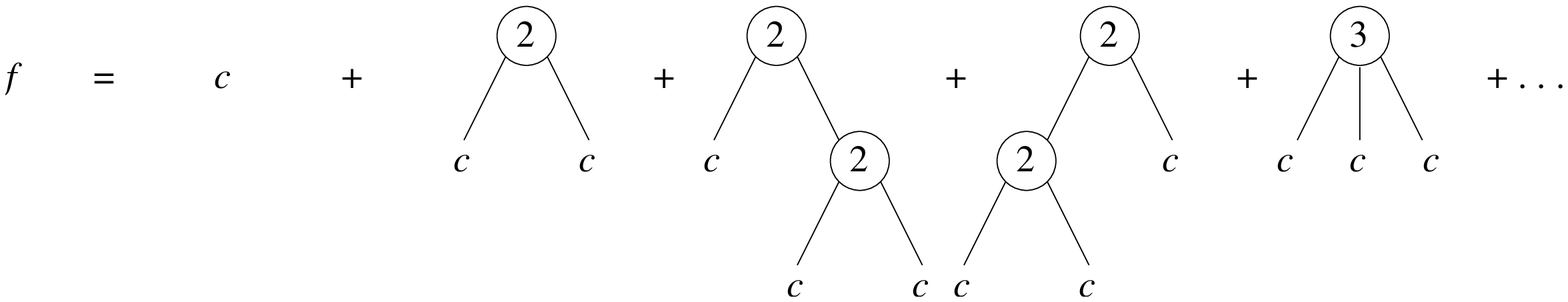}
\caption{\label{fig-premiersF}The terms $f_0$, $f_1$, $f_2$ expressed as
a sum of ordered trees.}
\end{center}
\end{figure}

\begin{figure}[ht]
\begin{center}
\leavevmode
\epsfxsize =10cm
\epsffile{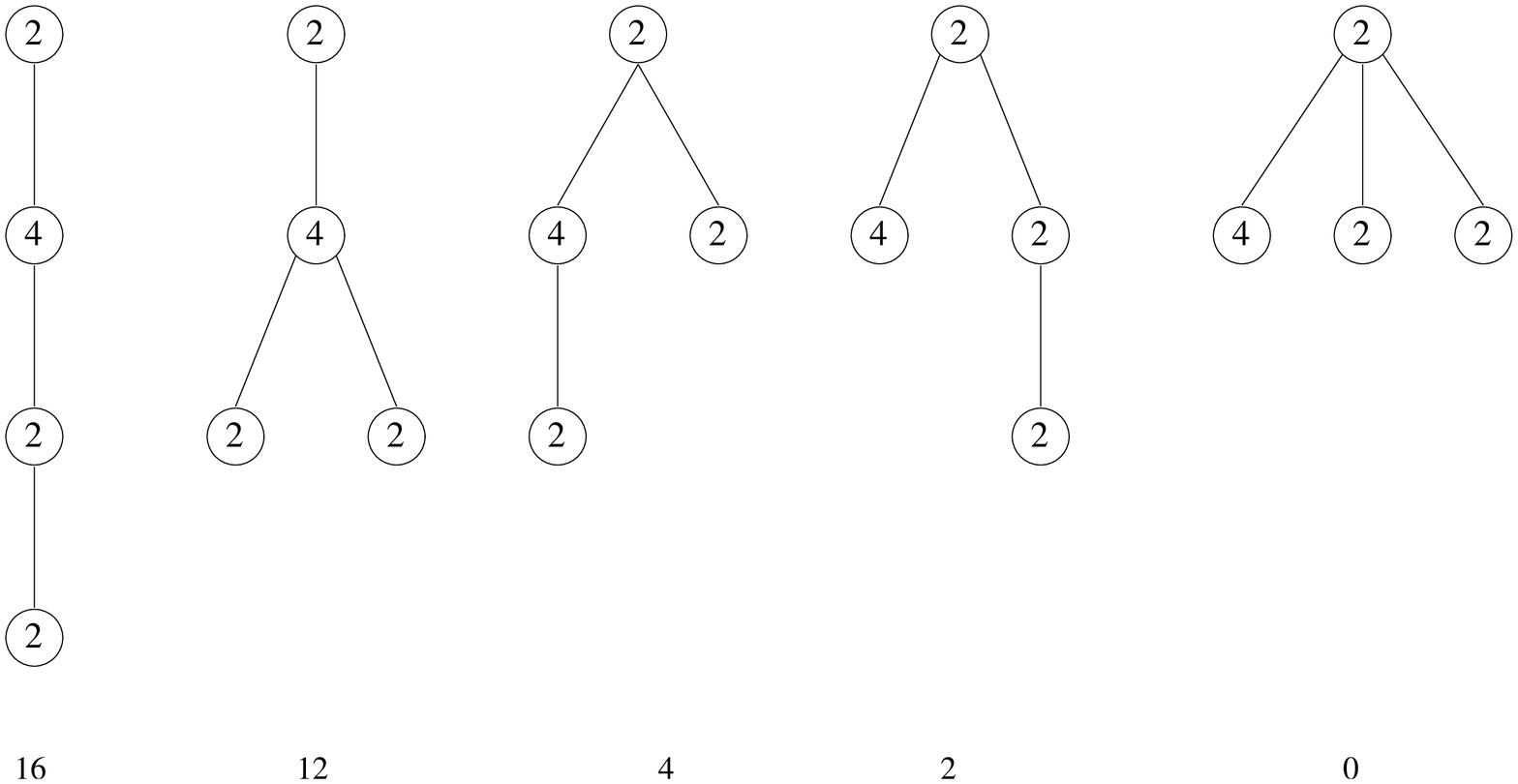}
\caption{\label{fig-ex34}}
\end{center}
\end{figure}

One can give for $h(A)$ a combinatorial interpretation analogous
to (\ref{gensol}).
Starting from the generalized inversion problem
\begin{equation}
f^{-1} = S_0 + S_1f + S_2f^2 + \ldots,
\end{equation}
we recast it in the form
\begin{equation}
\label{feq2}
f = c + d_1f^2 + d_2f^3 + \ldots,
\end{equation}
setting $c=S_0^{-1}$ and $d_n = - S_0^{-1} S_n$.
Solving recursively for $f_0$, $f_1$, $\ldots$, we find
\begin{equation}
f_0 =c, \quad f_1 = d_1cc, \quad
f_2 = d_1cd_1cc + d_1d_1ccc + d_2ccc\,,
\end{equation}
and we can now interpret each $d_i$ as the symbol of 
an $(i+1)$-ary operation in Polish notation.
Then,  $f_n$ is the sum of Polish codes of ordered trees with no
vertex of arity $1$ on $2n+1$ vertices (or Schr\"oder bracketings of the
words $c^{n+1}$) as one can check on Figure~\ref{fig-premiersF}.
From this, we can easily recover Formula~(2.21) of~\cite{BFK}. This 
amounts to set $c=1$, that is, solving
\begin{equation}
h = 1 + d_1h^2 + d_2h^3 + \ldots,
\end{equation}
as
\begin{equation}
h_n = \sum_{I\models n} \lambda_I d^I.
\end{equation}
We proceed as in the previous section.
Given the skeleton $S$ of a tree $T$, define its \emph{$1$-composition}
$I_1(S)$ as the sequence of values of the labels of the vertices of $S$ minus
$1$ in prefix order.
Then, thanks to Equation~(\ref{arb-squelette}), the number of trees with
skeleton $S$ is
\begin{equation}
\label{arb-squelette2}
\prod_{k=1}^p \binom{i_k+1}{a_k} 
\end{equation} 
where $I_1(S)=(i_1,\ldots,i_p)$ and $a_k$ is the arity of the $k$-th vertex
of $S$, numbered in prefix order.
For example, there are $34$ trees whose skeleton have $(1,3,1,1)$ 
as associated $1$-composition as one can check on Figure~\ref{fig-ex34}.

Let $I=(i_1,\ldots,i_p)$ be a composition of $n$.
The coefficient of $S^I$ in $h_n$ is equal to the number of ordered trees on
$2n+1$ vertices whose sequence of non-zero arities minus one in the prefix
reading is $I$. As before, the skeletons of these trees are
the ordered trees on $p$ vertices labeled by one plus the elements of $I$ in
prefix order.
The sequences of arities of the skeletons are the same as before so that
\begin{equation}
\label{lambda-res}
\lambda_I = \sum_{(a_1,\ldots,a_{p-1})}
\prod_{k=1}^{p-1} \binom{i_k+1}{a_k}
\end{equation}
where the sum is taken over the set of $a_k$ such that
$a_1+\cdots+a_j\geq j$ for all $j$ and $a_1+\cdots+a_{p-1} = p-1$.
This is Formula~(2.23) of~\cite{BFK}.
In this presentation, it is clear that the sum is over a Catalan set,
here the ordered trees.

\section{Explicit expressions in various bases}

\subsection{}
We shall now compute the coefficients of the expansions of
$g$ (or $h$ as well) in the bases $(R_I)$ and $(\Lambda^I)$
of $\Sym$.

The expansion on ribbons can be given for the $q$-analogs.
Let $\QRP(I)$ be the set of parking quasi-ribbons of shape $I$
(see \cite{NT1}), and let 
\begin{equation}
c_I(q)=\sum_{\a\in\QRP(I)}q^{\|\a\|}\,. 
\end{equation}
Then, since two words with the same evaluation
are hypoplactically equivalent iff the inverses of their standardized
have the same descents,
\begin{equation}
\ch_q(\PF_n)=G_n(q;A)=\sum_{I\vDash n}c_I(q) R_I(A)\,.
\end{equation}
For example,
\begin{eqnarray}
G_3(q;A) &=& S^{3}+ (q+q^2)S^{21}+{q}^{2}S^{12}+{q}^{3}S^{111}\nonumber\\
&=& (1+q+2q^2+q^3)R_3+(q+q^2+q^3)R_{21}+(q^2+q^3)R_{12}+q^3R_{111}
\end{eqnarray}
For $q=1$, this expansion presents a remarkable symmetry. The expansion
on elementary functions is given by the same formula as the expansion
on ribbons, up to sign and conjugation of the compositions:
\begin{equation}
g_n(A)=\sum_{I\vDash n}(-1)^{n-l(I)}c_{I^{\sim}}\Lambda^I\,.
\end{equation}
For example,
\begin{eqnarray}
g_3(A)&=&\Lambda^3-3\Lambda^{21}-2\Lambda^{12}+5\Lambda^{111}\,,\\
g_4(A)&=&-\Lambda^4+4\Lambda^{31}+3\Lambda^{22}+2\Lambda^{13}
  -9\Lambda^{211}-7\Lambda^{121}-5\Lambda^{112}+14\Lambda^{1111}\,.
\end{eqnarray}

This symmetry is equivalent to the invariance of $g$ under the
linear involution of $\Sym$ defined by 
\begin{equation}
\nu:\ S^I \longmapsto S^{I^\sim}\,,
\end{equation}
as one can check on Equation~(\ref{g01234}).
Indeed,
\begin{lemma}
On the ribbon basis, $\nu$ is given by
\begin{equation}
\nu(R_I)=(-1)^{l(I)-1}\Lambda^{I^\sim}\,.
\end{equation}
\end{lemma}

\Proof The image of the Cauchy kernel $\sigma_1(XA)$ by $\nu$
is
\begin{equation}
\begin{split}
\sum_I F_I \nu(R_I) &= \sum_I M_I \nu(S^I) = \sum_I M_I S^{I^\sim}=\sum_I
M_{I^\sim}S^I\\
&= \sum_I\sum_{J\le I}M_{I^\sim}R_J=\sum_I M_{I^\sim}
\sum_{J^\sim\ge I^\sim}R_{J^\sim}\\
&=\sum_I M_{I^\sim}\sum_{J^\sim\ge I^\sim}\sum_{K\le J}(-1)^{l(J)-l(K)}\Lambda^K\\
&=\sum_K\left( \sum_{J\ge K}(-1)^{l(J)}\sum_{I^\sim\le J^\sim}M_{I^\sim}
\right)\Lambda^K\\
&=\sum_K\left(\sum_I M_{I^\sim}\sum_{K\le J\le I}(-1)^{l(J)}\right)
(-1)^{l(K)}\Lambda^K\\
&=\sum_K\left(\sum_{I\ge K}M_{I^\sim}(-1)^{l(K)}\right)\Lambda^K
=\sum_K F_{K^\sim}(-1)^{l(K)}\Lambda^K\\
&=\sum_K F_K(-1)^{l(K)-1}\Lambda^{K^\sim}\,.
\end{split}
\end{equation}
\qed

\subsection{An involution}

Actually, the $\nu$-invariance of $g$ follows from a stronger property.
As we have seen, the solution $f$ of the general inversion problem
(\ref{genfeq})
\begin{equation}
f=S^0 + S^{10} + S^{200}+S^{110}+
S^{1110}+S^{1200}+S^{2010}+S^{2100}+S^{3000}+\cdots
\end{equation}
can be interpreted as the formal sum of all nondecreasing parking functions.
We will now prove that there exists a canonical involution $\iota$ on these which
satisfies
\begin{equation}
\evt(\iota(\pi))=\evt(\pi)^\sim\,.
\end{equation}
To simplify the presentation, we shall identify a nondecreasing parking
function with its evaluation. 
More precisely,  define a \emph{generalized composition} as a composition
where zeros are allowed. The composition obtained by removing all zeros is
called the \emph{corresponding composition}.
A generalized composition $I$ of $n$ is of \emph{parking type} iff it is
of length $n+1$ and $i_1+\ldots+i_k\geq k$ for all $k$ in $[1,n]$. In other
words, the set of generalized compositions of parking type is the set of
evaluations of parking functions with an appropriate number of trailing zeros.

\begin{figure}[ht]
\begin{center}
\leavevmode
\epsfxsize =5cm
\epsffile{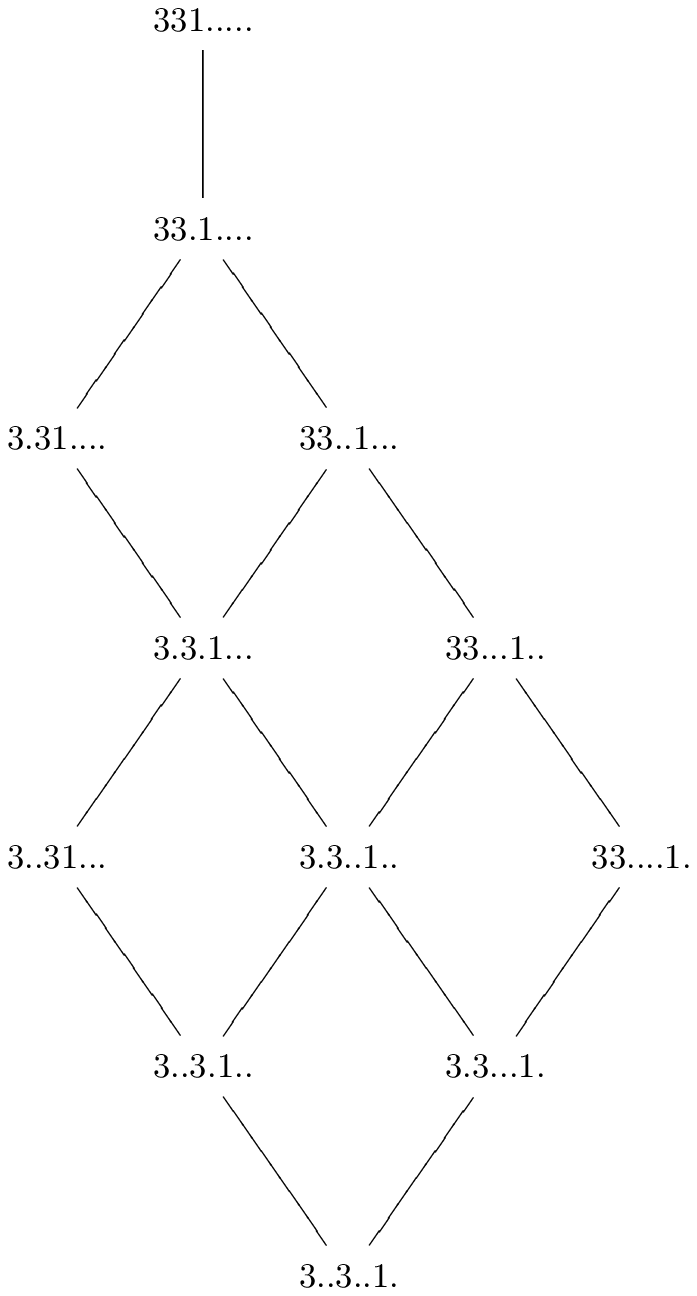}
\hskip3cm
\epsfxsize =5cm
\epsffile{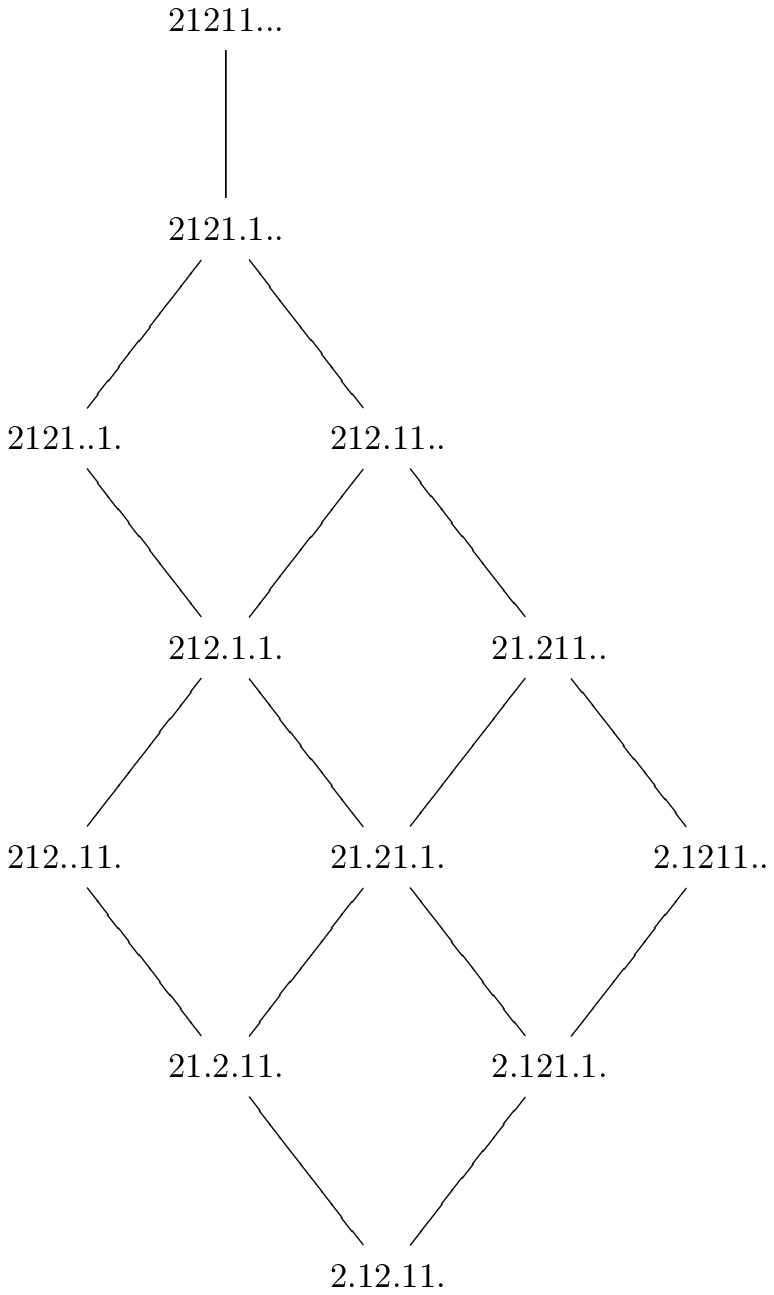}
\end{center}
\caption{\label{G331}The graphs $\Gamma_{331}$ and $\Gamma_{21211}$  ($0$s are
represented by dots).}
 \end{figure}

Before describing the involution on generalized compositions of parking type,
we need  some more structure on the set of elements having the same packed
evaluation. For each composition $I$ of $n$, build a directed graph $\Gamma_I$
with vertex set given by generalized compositions of parking type with
corresponding composition $I$
and an arrow $J \longrightarrow J'$ iff $J'$ is obtained from $J$ by
exchanging two consecutive parts of $J$, $j_i$ and $j_{i+1}$ so that $j_i$ or
$j_{i+1}$ is $0$, an operation reminescent of Hivert's quasi-symmetrizing
action~\cite{Hiv}.
For example, $\Gamma_{331}$ and $\Gamma_{21211}$ are given on
Figure~\ref{G331}.
$\Gamma_I$ can be seen as an initial interval of a permutohedron:
consider the word $K=(0^{n+1-l(I)})$ and the shuffle $S=I\shuffle K$. To these
elements corresponds naturally an element of the shuffle
$S'=(123\cdots l(I)) \shuffle (l(I)+1\cdots n+1)$.
Then if one restricts to the elements of $S$ that are the evaluation of a
nondecreasing parking function, for any such element $s$, all the permutations
smaller than the corresponding element $s'$ in $S'$ 
correspond to  evaluations of nondecreasing parking functions: indeed,
this  means that if $J$ is of parking type, all generalized compositions
obtained from $J$ by moving zeros to the right also are of parking type,
which is obvious. Now there is only one minimal element, the concatenation of
$I$ and $K$, and only one maximal element, the evaluation where any non-zero
entry $i$, except for the last one, is followed by exactly $i-1$ zeros before
the next non-zero entry: no successor of this element is the evaluation of a
parking function and all other elements have at least one successor of this
type.

We are now in a position to describe the involution on generalized
compositions of parking type which induces the conjugation on the underlying
compositions.

\begin{algorithm}
Let $J$ be a generalized composition of parking type.
\begin{itemize}
\item Let $J'$ be the tuple obtained by reading $J$ from right to left.
\item Compute the conjugate $C$ of the corresponding composition of $J$, 
\item fill the zero slots of $J'$ by the parts of $C$, 
\item replace by $0$ the nonzero parts of $J'$.
\item This is the output, $\iota(J)$.
\end{itemize}
\end{algorithm}

For example, if $J=(2,1,1,0,1,2,0,2,0,0)$, then
$J'=(0,0,2,0,2,1,0,1,1,2)$, the corresponding
composition being $(2,1,1,1,2,2)$ and its conjugate $(1,2,5,1)$. 
We insert the four parts of this last composition into the zero slots of $J'$,
putting zeros at the other places and get $(1,2,0,5,0,0,1,0,0,0)$.

\begin{lemma}
The previous algorithm is an involution on generalized compositions of parking
type, sending maximal elements of graphs to maximal elements.
\end{lemma}

\Proof
The algorithm is an involution since the conjugation of compositions is one, so we
only have to prove that the output is of parking type if the input is.
By construction of $\Gamma_I$, it is sufficient to prove that the image of the
bottom element of $\Gamma_I$ is of parking type. Thanks to its
characterization, it is obvious that this bottom element is sent to the bottom
element of $\Gamma_{I^\sim}$ by our involution.
\qed

\begin{theorem}
The graphs $\Gamma_I$ and $\Gamma_{I^\sim}$ associated to mutually
conjugate compositions of $n$ are isomorphic.
Moreover, if one labels the edges by $i$ when one exchanges the letters
in positions $i$ and $i+1$, then the labels of the edges are exchanged by the
involution $i\leftrightarrow n+1-i$.
\end{theorem}

\Proof
The graph $\Gamma_I$ corresponds to a part of the shuffle
$I\shuffle 0^{n+1-l(I)}$ whereas the graph $\Gamma_{I^\sim}$ corresponds to a
part of the shuffle $I^\sim \shuffle 0^{n+1-l(I^\sim)}$. It is known that 
$l(I)+l(I^\sim)=n+1$ so that both graphs correspond to parts of a shuffle of
an element of length $l(I)$ with an element of length $l(I^\sim)$.
Moreover, given the definition of the edges of both graphs, an edge labelled
$i$ between $P$ and $P'$ proves that there is an edge labelled $n+1-i$
between $\iota(P)$ and $\iota(P')$. So both graphs are isomorphic.
\qed

For example, the two graphs on Figure~\ref{G331} corresponding to $331$ and
$331^\sim$ are indeed isomorphic.

\subsection{A representation theoretical interpretation}

In fact, nondecreasing parking functions form a sub-semigroup of the semigroup
of all endofunctions of $[n]$. Its representation theory has been investigated
by Hivert and Thi\'ery \cite{HT}, and it follows from their work that the
graphs $\Gamma_I$ (now seen on nondecreasing parking functions instead of 
generalized  compositions of parking type) encode the indecomposable
projective modules $P_I$ of the semigroup algebra
${\mathcal C}_n=\C[\NDPF_n]$.
Indeed, these modules are parametrized by compositions of $n$, and 
each $P_I$ has a basis $(b_\pi)_{\pi\in\NDPF_I}$, such that if one denotes
by $e_i$ the generator mapping $i+1$ to $i$ and leaving invariant all other
$j$, $e_i\circ b_\pi=b_{\pi'}$ iff $\pi \overset{i}{\longrightarrow}\pi'$ and  
$e_i\circ b_\pi=0$ otherwise.
Thus, on the one hand, the coefficients $\delta_I$ of the expansion
\begin{equation} \label{dISI}
g_n=\chb(\PF_n)=\sum_{I\vDash n} \delta_I S^I
\end{equation}
are the dimensions $\delta_I=\dim\, P_I$ of the indecomposable
projective modules of ${\mathcal C}_n$. On the other hand, the noncommutative
symmetric functions $S^I$ are the characteristics of the permutational modules
of $H_n(0)$, which are projective, but decomposable for $I\not=(n)$.
As also shown in \cite{HT}, these permutational modules are in fact the
indecomposable projective modules for a larger algebra, the {\em
Hecke-symmetric} algebra $\HS_n$. One can check that the right action
of $\HS_n$ on $\C\PF_n$ and the left action of ${\mathcal C}_n$
(by composition $\pi\circ\a$) commute with each other, so that
the expression (\ref{dISI}) of $\chb(\PF_n)$ reflects the decomposition
of $\C\PF_n$ as a $({\mathcal C}_n,\HS_n)$-bimodule.

The coefficients $(\lambda_I)$ of the ribbon expansion can be similarly
interpreted as the dimensions of the projective modules of the commutant of
$H_n(0)$ in $\C\PF_n$, an algebra ${\mathcal D}_n$ having as dimension the
Schr\"oder number $s_n$ and containing ${\mathcal C}_n$.

\section{Noncommutative Abel identities}

\subsection{}
Abel's generalization of the binomial identity can be stated as
\begin{equation}
p_n(x+y)=\sum_{k=0}^n\binom{n}{k}p_k(x)p_{n-k}(y)\,,
\end{equation}
that is, the Abel polynomials
\begin{equation}
p_n(x):=x(x+n)^{n-1}
\end{equation}
form a sequence of binomial type.

Let $\E$ be the exponential alphabet. If we set $g_n=\chb(\PF_n)$,
we have $g_n(t\E)=(n+1)^{n-1} t^n/n!$, and it follows
from  (\ref{lem1}) that
\begin{equation}
g^r=\sum_{n\ge 0} \chb(\PF_n^{(r)}) =: g^{(r)}
\end{equation}
where 
\begin{equation}
g_n^{(r)}(t\E)= \frac{t^n}{ n!}|\PF_n^{(r)}|=\frac{t^n}{n!} r(r+n)^{n-1}
\end{equation}
and $\PF_n^{(r)}$ is the set of words $\a$ whose nondecreasing
rearrangement satisfies $\a^\uparrow_i\le i+r-1$.
This is because of the self-evident generalization of (\ref{lem1})
\begin{equation}
\PF_n^{(r)} = \bigsqcup_{k=0}^n \PF_k^{(j)}\, \shuffle \PF_{n-k}^{(r-j)}\,,
\end{equation}
for all $j$ in $[1,r-1]$, which implies in particular
\begin{equation}
\chb (\PF_n^{(r)}) = \sum_{k=0}^n \chb(\PF_k) \chb(\PF^{(r-1)}_{n-k}).
\end{equation}

Hence, the $t\E$-specialization of $g^x$ is the exponential generating
function of Abel's polynomials, and Abel's identity amounts to the obvious
equality $g^xg^y=g^{x+y}$.

We can therefore define the noncommutative Abel polynomial $P_n(x;A)$ as the
term of degree $n$ in $g(A)^x$. It can be computed directly 
using the  binomial expansion of $g(A)^x=(1+U)^x= \sum_n\binom{x}{n}U^n$.
For example,
\begin{equation}
\begin{split}
P_1(x;A) =& xS^1\,,\\
P_2(x;A) =& xS^{2} + \frac{x(x+1)}{2}S^{11}\,,\\
P_3(x;A) =& xS^{3} + \frac{x(x+3)}{2}S^{21} + \frac{x(x+1)}{2}S^{12}
+ \frac{x(x+1)(x+2)}{6}S^{111}\,, \\
P_4(x;A) =& xS^{4} + \frac{x(x+5)}{2}S^{31} + \frac{x(x+3)}{2}S^{22}
           + \frac{x(x+1)}{2}S^{13} \\
&+ \frac{x^3+6x^2+11x}{6}S^{211}
           + \frac{x^3+6x^2+5x}{6}S^{121}\\
&           + \frac{x^3+3x^2+2x}{6}S^{112}
           + \frac{x(x+1)(x+2)(x+3)}{6}S^{1111}.
\end{split}
\end{equation}

In particular, one has
\begin{equation}
\chb (\PF_n^{(r)})=P_n(r;A)\,.
\end{equation}
But this characteristic can also be computed directly. Indeed, since
$\PF_n^{(r)}$ is a permutational module, we have
\begin{equation}
\chb (\PF_n^{(r)})=\sum_{I\vDash n}\alpha_IS^I\,,
\end{equation}
where $\alpha_I$ is the number of nondecreasing words $\a\in\PF_n^{(r)}$
with packed evaluation $I$. These elements can be classified according
to their parkized $\b=\Park(\a)$ (see \cite{NT1}), which is an
ordinary  nondecreasing parking function. The cardinality $\alpha_\b$ 
of such a class is a binomial coefficient. To see this, let
\begin{equation}
\b=\b_1\bullet\b_2\bullet\cdots\bullet\b_m
\end{equation}
be the maximal factorization of $\b$ into connected nondecreasing
parking functions ($\bullet$ denoting shifted concatenation, see \cite{NT1}).
The nondecreasing $\a\in\PF_n^{(r)}$ such that $\Park(\a)=\b$ are
obtained by shifting each factor $\b_i$ of an amount $k_i$,
such that $k_1+\cdots+k_m\le r$. Thus,
\begin{equation}\label{alphab}
\alpha_\b=\binom{r+m-1}{m}\,.
\end{equation} 
Set $c(\b)=m$. Formula (\ref{alphab}) being valid for all positive
integers $r$, we have in general
\begin{equation}\label{ncab}
P_n(x;A)=\sum_{I\vDash n}\left(
\sum_{\b\in\NDPF_n;\ \evt(\b)=I}\binom{x+c(\b)-1}{c(\b)}\right)S^I\,.
\end{equation}
For example, there are three nondecreasing parking functions
with packed evaluation $(211)$: $1123$, $1124=112\bullet 1$
and $1134=11\bullet 1\bullet 1$, so that
\begin{equation}
\alpha_{211}=\binom{x}{1}+\binom{x+1}{2}+\binom{x+2}{3}
=\frac{x^3+6x^2+11x}{6}\,.
\end{equation}
Similarly, the coefficient of $S^{31}$ in $P_3$ is
\begin{equation}
\alpha_{31}=\binom{x}{1}+\binom{x}{1}+\binom{x+1}{2}=\frac{x(x+5)}{2}\,,
\end{equation}
corresponding to $1112$, $1113$ and $1114=111\bullet 1$.

\subsection{}
By construction, the specialization $A=\E$ gives back
the Abel polynomials.
As usual, the specialization $A=1$ is also interesting.
Let $a(n,m)$ be the Catalan triangle \cite[A009766]{Slo}. That is,
\begin{equation}
a(n,m)= \binom{n+m}{n}\frac{n-m+1}{n+1}\,
\end{equation}
and 
\begin{equation}
\sum_{n\ge 0}\left(\sum_{m=0}^n a(n,m)t^m\right)z^n = \frac{C(tz)}{1-zC(tz)}\,,
\end{equation}
where 
\begin{equation}
C(z)=\frac{1-\sqrt{1-4z}}{2z}
\end{equation}
is the generating series of the Catalan numbers.
We need the reverted and shifted triangle
\begin{equation}
c(n,k)=a(n-1,n-k-1) \ (n\ge 1)\,,\ \ c(0,0)=1\,,
\end{equation}
whose generating series is
\begin{equation}\label{trcat}
\frac{1}{1-tzC(z)}= 1+tz +(t^2+t)z^2+(t^3+2t^2+2t)z^3+\cdots
=\sum_{n\ge 0}\left(\sum_{m=0}^n c(n,m)t^m\right)z^n\,.
\end{equation}
Let 
\begin{equation}
S_n(x)=\frac{x(x+1)\cdots(x+n-1)}{n!}\,,
\end{equation}
that is, $S_n(x)$ is the coefficient of $t^n$ in $(1-t)^{-x}$
(which can be interpreted as $\sigma_t(x)$ for $x$ a binomial
element, whence the choice of notation).

We can now state:
\begin{theorem}
The specialization $A=1$ of the noncommutative Abel
polynomials $P_n(x;A)$ is given by
\begin{equation}\label{Pncnk}
P_n(x;1)=\sum_{k=1}^n c(n,k)S_k(x)\,.
\end{equation}
Moreover, their generating series is
\begin{equation}\label{gfPnA1}
\sum_{n\ge 0}P_n(x;1)z^n = C(z)^x\,.
\end{equation}
\end{theorem}

\Proof Equation (\ref{gfPnA1}) is clear if one rewrites
the quadratic equation for $C(z)$ as
\begin{equation}\label{feqC}
C(z)= \frac{1}{1-zC(z)}=\sum_{n\ge 0}S_n(1)[zC(z)]^n\,.
\end{equation}
Equation (\ref{Pncnk}) follows from (\ref{ncab}), since
(\ref{trcat}) shows that $c(n,k)$ is the number of
nondecreasing parking functions of length $n$ such that $c(\b)=k$.
It can also be proved analytically. 
The generating series of the right-hand sides of (\ref{Pncnk})
can be written as a contour integral, over a circle
$\gamma=\{|w|=\varepsilon <1\}$
\begin{equation}
\frac{1}{2\pi i}\oint_\gamma \frac{(1-w)^{-x}}{w-zC(z)}dw
= \frac{1}{2\pi i}\oint_\gamma \frac{f(w)}{w-a}dw
\end{equation}
where $a=zC(z)$ and $f(w)=(1-w)^{-x}$. For $|z|$ small enough,
$a$ is inside $\gamma$, and by Cauchy's theorem, the right-hand side
is
\begin{equation}
f(a)= (1-zC(z))^{-x} = C(z)^x
\end{equation}
according to (\ref{feqC}). 
\qed

The coefficients of the $P_n(x;1)$ build up the triangle
\cite[A038455]{Slo}.
In fact, $C(z)={\mathcal B}_2(z)$,
where
\begin{equation}
{\mathcal B}_t(z)= \sum_{n\ge 0}(tn)^{\underline{n-1}}\cdot \frac{z^n}{n!}
\end{equation}
is Lambert's generalized binomial series (see \cite{GKP}, (5.68) p. 200).
According to \cite{GKP}, (5.70), we have finally the closed expression
\begin{equation}
P_n(x;1)=\binom{x+2n}{n}\frac{x}{x+2n}\,.
\end{equation}

\section{$(k,l)$-Parking functions}

There is a general notion of parking functions associated with
a sequence ${\bf u}=(u_n)_{n\ge 1}$ of positive integers: these
are the words $\a$ such that $(\a^\uparrow)_i\le u_i$.
In general, their enumeration can be obtained only in terms of
Gon\u carov polynomials \cite{KY}. In the particular case where
${\bf u}$ is an arithmetic progression, it is possible to obtain closed
formulas, of which we shall now give the noncommutative analogs.

Let 
\begin{equation}
\PF^{(k,l)}_n=\{\a\in[l+(n-1)k]^n\,|\, \a^\uparrow_i\le l+(i-1)k\}\,.
\end{equation}
Stanley and Pitman \cite{PS} have shown that
\begin{equation}
|\PF^{(k,l)}_n|=l(l+kn)^{n-1}\,.
\end{equation}
As above, this can be extended to the calculation of the $0$-Hecke
characteristic. The argument used for  (\ref{lem1}) proves as well
(cf. \cite{KY})
\begin{equation}\label{lem2}
[N]^n=\bigsqcup_{j=0}^n \PF^{(k,l)}_j \shuffle [jk+l+1,N]^{n-k}\,.
\end{equation}
Taking $N=nk+r$, we obtain for the characteristic
\begin{equation}
\sum_{n\ge 0}S_n((nk+r)A)=\left(\sum_{n\ge0}\chb\, \PF^{(k,l)}_n\right)
\left(\sum_{n\ge0}S_n((nk+r-l)A)\right)\,.
\end{equation}
Setting
\begin{equation}
F(r,k)=\sum_{n\ge 0}S_n((nk+r)A)
\end{equation}
we have finally
\begin{proposition}
The  noncommutative characteristic of the permutational $H_n(0)$-module
on $\PF^{(k,l)}_n$ is equal to $g^{(k,l)}_n(A)$,
the term of degree $n$ in $g^{(k,l)}(A)=F(r,k)F(r-l,k)^{-1}$, which is independent of $r$.
\end{proposition}

\begin{example}{\rm Taking $A=t\E$ and $r=l$, we recover the enumeration
\begin{equation}
\sum_{n\ge 0}|\PF^{(k,l)}_n|\frac{t^n}{n!}=
\frac{\sum_{n\ge 0}\frac{t^n}{n!}(nk+l)^n}
     {\sum_{n\ge 0}\frac{t^n}{n!}(nk)^n}
=\sum_{n\ge 0}\frac{t^n}{n!}l(nk+l)^{n-1}
\end{equation}
the last equality following from Abel's identity (the middle term
is $g(tk\E)^{l/k}$). Note that this can also be expressed in terms
of the generalized exponential series of \cite{GKP}
\begin{equation}
{\mathcal E}_\alpha(z)=\sum_{n\ge 0}(n\alpha+1)^{n-1}\frac{z^n}{n!}\,,
\end{equation}
that is,
\begin{equation}
\sum_{n\ge 0}|\PF^{(k,l)}_n|\frac{t^n}{n!}=
{\mathcal E}_{\frac{k}{l}}(lt)\,.
\end{equation}
If we set 
\begin{equation}
{\mathcal E}(z)= g(z\E)\,,
\end{equation}
we see that
\begin{equation}
{\mathcal E}_{\frac{k}{l}}\left(\frac{l}{k}t\right)={\mathcal
E}\left(k\frac{t}{k}\right)^{\frac{l}{k}}\,.
\end{equation}
Hence, for any $\alpha$
\begin{equation}\label{Ealpha}
{\mathcal E}_\alpha(t)={\mathcal E}(\alpha t)^{\frac1\alpha}
\end{equation}
since this is true for $\alpha$ rational. This equality implies most
of the interesting properties of ${\mathcal E}_\alpha(z)$. Indeed,
let us write down explicitly the functional equation for $g(z\E)$
\begin{equation}
g(z\E)=\sum_{n\ge 0}S_n( z\E)g( z\E)^n=
\sum_{n\ge 0}\frac{ z^n}{n!}g( z\E)^n=e^{ z g( z\E)}\,.
\end{equation}
We see that 
${\mathcal E}(z)$ is Eisenstein's function, defined by
\begin{equation}
{\mathcal E}(z)=e^{z{\mathcal E}(z)}
\end{equation}
(see \cite{GKP}, (5.68) p. 200). Now, (\ref{Ealpha})
implies immediately identities like \cite{GKP}, (5.69)
\begin{equation}
{\mathcal E}_\alpha(z)^{-\alpha}\ln {\mathcal E}_\alpha(z)=z
\end{equation}
for instance.

}
\end{example}

\begin{example}{\rm The specialization $A=1$ gives
\begin{equation}
g^{(k,l)}_n(1)= \binom{nk+l+n-1}{n}\frac{l}{nk+l}\,.
\end{equation}
This is defined for negative values of $k$ and $l$ as well,
and we have
\begin{equation}
\sum_{n\ge 0}g^{(-k,-l)}_n(1)(-t)^n=
\sum_{n\ge 0}\binom{nk+l}{n}\frac{l}{nk+l}t^n ={\mathcal B}_k(t)^l
\end{equation}
according to \cite{GKP}, (5.70). Again, this can be generalized.
Recall that ${\mathcal B}_\alpha(z)$
is defined by
\begin{equation}
{\mathcal B}_\alpha(z)=1+z\sum_{n\ge 0}\binom{(n+1)\alpha}{n}\frac{z^n}{n+1}
=1+z\sum_{n\ge 0}\frac{\Lambda^n((n+1)\alpha)}{n+1}z^n\,, 
\end{equation}
(using the $\lambda$-ring notation), so that
\begin{equation}
{\mathcal B}_{-\alpha}(-z)=
1-z\sum_{n\ge 0}\frac{S^n((n+1)\alpha)}{n+1}z^n = 1-zg(z\alpha)\,,
\end{equation}
where $g(z\alpha)$ denotes the specialization $A=z\alpha$ ($\alpha$
binomial and $z$ of rank 1) of $g(A)$. Hence, $g(z\alpha)$ satisfies
the functional equation
\begin{equation}
g(z\alpha)=\sum_{n\ge 0}S_n(z\alpha) g(z\alpha)^n=(1-zg(z\alpha))^{-\alpha}\,,
\end{equation}
so that
\begin{equation}\label{Bandg}
g(z\alpha)={\mathcal B}_{-\alpha}(-z)^{-\alpha}\,.
\end{equation}
Clearly, this implies
\begin{equation}
{\mathcal B}_{-\alpha}(-z)^{1+\alpha}-{\mathcal B}_{-\alpha}(-z)^{\alpha}=-z\,,
\end{equation}
which is the first equation of \cite[(5.69)]{GKP}.
Hence, the specialization $A=1$ explains the properties of the generalized
binomial
series ${\mathcal B}_\alpha(z)$, while $A=\E$ takes cares of the generalized
exponential series ${\mathcal E}_\alpha(z)$.
Note that (\ref{Bandg}) allows one to write
\begin{equation}
{\mathcal B}_\alpha(-z)=g(z\alpha)^{\frac{1}{\alpha}}
\end{equation}
(but {\it not} ${\mathcal B}_\alpha(z)=g(-z\alpha)^{\frac{1}{\alpha}}$,
since in the right-hand side, $g$ is interpreted as a $\lambda$-ring
operator!). 
}
\end{example}

The $q$-characteristic is obtained similarly:
\begin{proposition}
The $q$-characteristic of the permutational $H_n(0)$-module
on $\PF^{(k,l)}_n$ admits as generating series
\begin{equation}
\sum_{n\ge 0}x^n q^{-k\binom{n+1}{2}-n(nk+l)}\chb_q \PF^{(k,l)}_n
= F_{k }^{(r)}(x,q;A)F_{k }^{(r-l)}(q^lx,q;A)^{-1}
\end{equation} 
where
\begin{equation}
F_{k }^{(r)}(x,q;A)=\sum_{n\ge 0}x^n q^{-k\binom{n+1}{2}}S_n([nk+r]_qA)\,.
\end{equation}
\end{proposition}

\begin{example}{\rm Taking $A=1$, $k=3,l=2$ and $r=\infty$, we obtain
a $q$-analog of sequence \cite[A069271]{Slo} (1,2,9,52,340,2394,17710,...) 
for the $q$-enumeration of nondecreasing $\PF^{(3,2)}$'s
(cf. \cite{Yan}, Cor. 5.1.)
\begin{equation}
\begin{split}
f(t,q)&=1+ (q+1)q^{-3}t+ (q^5+2q^4+2q^3+2q^2+q+1)q^{-9}t^2\\
&+ (q^{12}\!+\!3q^{11}\!+\!5q^{10}\!+\!7q^9\!+\!7q^8\!+\!7q^7\!
+\!6q^6\!+\!5q^5\!+\!4q^4\!+\!3q^3\!+\!2q^2\!+\!q\!+\!1) q^{-18}t^3\\
&+O\left (t^4\right ).
\end{split}
\end{equation}
}
\end{example}

\section{Generalized inversion formulas}

\subsection{Some families of trees}
Let $b$ be an integer.
Consider the generalized inversion problem
\begin{equation}
\label{general-f}
f = c + d_1 f^{b+1} + d_2 f^{b+2} + \ldots
\end{equation}
which reduces to Equations~(\ref{genfeq}) and~(\ref{feq2}) for $b=0$ and
$1$.

The general solution $f_n$ is given by the Polish codes of ordered trees with
$(b+1)n+1$ leaves and no vertices of arity between $1$ and $b$.

We are interested in the special case
\begin{equation}
\label{general-g}
g = 1 + d_1 g^{b+1} + d_2 g^{b+2} + \ldots
\end{equation}
that is, we want to compute the coefficients $\delta^{(b)}_I$ of
\begin{equation}
g_n = \sum_{I\models n} \delta^{(b)}_I d^I.
\end{equation}

Given the skeleton of a tree such that each internal vertex is of arity at
least $b+1$, one can define its $b$-composition 
$I_b(S)$ as the sequence of values of the labels of the vertices of $S$ minus
$b$, in prefix order.
Thanks to Equation~(\ref{arb-squelette}), the number of trees with skeleton
$S$ is
\begin{equation}
\label{arb-squeletteg}
\prod_{k=1}^p \binom{i_k+b}{a_k} 
\end{equation} 
where $I_b(S)=(i_1,\ldots,i_p)$ and $a_k$ is the arity of the $k$-th vertex
of $S$ in prefix order, so that, as in Equations~(\ref{delta-res})
and~(\ref{lambda-res}):

\begin{equation}
\delta^{(b)}_I = \sum_{(a_1,\ldots,a_{p-1})}
\prod_{k=1}^{p-1} \binom{i_k+b}{a_k}
\end{equation}
where the sum is again taken over sequences $(a_1,\ldots,a_{p-1})$ such that
$a_1+\cdots+a_j\geq j$ for all $j$ and $a_1+\cdots+a_{p-1} = p-1$.
%

\subsection{Combinatorial triangles}

Let now $\gamma^{(b)}_{p,n}$ be
\begin{equation}
\gamma^{(b)}_{p,n} := \sum_{I\models n ; i_1=p} \delta^{(b)}_I.
\end{equation}
This amounts to enumerate the trees by  arity of the root.
The triangles $(\gamma^{(b)}_{p,n})$ include some classical triangles of the
combinatorial literature:
for $n=0$, one recovers the Catalan triangle (sequence A033184 of~\cite{Slo}),
for $n=1$, one recovers the Schr\"oder triangle (sequence A091370
of~\cite{Slo}). Their first terms are given on Figure~\ref{catschro}.

\begin{figure}[ht]
$$
\begin{array}{cccccccc}
  1 &     &    &    &    &   &   \\
  1 &   1 &    &    &    &   &   \\
  2 &   2 &  1 &    &    &   &   \\
  5 &   5 &  3 &  1 &    &   &   \\
 14 &  14 &  9 &  4 &  1 &   &   \\
 42 &  42 & 28 & 14 &  5 & 1 &   \\
132 & 132 & 90 & 48 & 20 & 6 & 1 \\
\end{array}
\qquad
\begin{array}{cccccccc}
   1 &      &     &     &    &   &  \\
   2 &    1 &     &     &    &   &  \\
   7 &    3 &   1 &     &    &   &  \\
  28 &   12 &   4 &   1 &    &   &  \\
 121 &   52 &  18 &   5 &  1 &   &  \\
 550 &  237 &  84 &  25 &  6 & 1 &  \\
2591 & 1119 & 403 & 125 & 33 & 7 & 1\\
\end{array}
$$
\caption{\label{catschro}The Catalan and Schr\"oder triangles ($b=0$ and
$b=1$).}
\end{figure}

The  triangles for $b=2$ and $b=3$ are given on  Figure~\ref{b=23}. Note that
although they are not (yet) referenced in~\cite{Slo}, the row sums of the
case $b=2$ yields  Sequence A108447, with a quite different interpretation.

\begin{figure}[ht]
$$
\begin{array}{cccccccc}
    1 &      &     &     &    &   &  \\
    3 &    1 &     &     &    &   &  \\
   15 &    4 &   1 &     &    &   &  \\
   85 &   22 &   5 &   1 &    &   &  \\
  519 &  132 &  30 &   6 &  1 &   &  \\
 3330 &  837 & 190 &  39 &  7 & 1 &  \\
22135 & 5516 &1250 & 260 & 49 & 8 & 1\\
\end{array}
\qquad
\begin{array}{cccccccc}
     1 &       &      &     &    &   &  \\
     4 &     1 &      &     &    &   &  \\
    26 &     5 &    1 &     &    &   &  \\
   192 &    35 &    6 &   1 &    &   &  \\
  1531 &   270 &   45 &   7 &  1 &   &  \\
 12848 &  2215 &  362 &  56 &  8 & 1 &  \\
111818 & 18961 & 3054 & 469 & 68 & 9 & 1\\
\end{array}
$$
\caption{\label{b=23}Triangles obtained for $b=2$ and $b=3$.}
\end{figure}

One can also choose negative values for $b$ even if this has no
direct interpretation in terms of trees. With $b=-1$, the equation becomes
\begin{equation}
\label{f0}
f = c + d_1 + d_2f + d_3f^2 + \ldots
\end{equation}
and one recovers up to sign the Motzkin triangle (sequence A091836
of~\cite{Slo}, see Figure~\ref{b=-1-2}) splitting up the Motzkin numbers
(sequence A001006 of~\cite{Slo}) when putting $c=1$ and considering $d_i$ as a
$(i\!-\!1)$-ary operation. Recall that Motzkin paths are the paths from
$(0,0)$ to $(n,0)$, with three kinds of steps $(1,0)$, $(1,1)$, and $(1,-1)$,
that never go below the horizontal axis.

\begin{figure}[ht]
$$
\begin{array}{cccccccc}
    1 &      &     &     &    &   &  \\
    1 &      &     &     &    &   &  \\
    1 &    1 &     &     &    &   &  \\
    1 &    2 &   1 &     &    &   &  \\
    2 &    3 &   3 &   1 &    &   &  \\
    4 &    6 &   6 &   4 &  1 &   &  \\
    9 &   13 &  13 &  10 &  5 & 1 &  \\
   21 &   30 &  30 &  24 & 15 & 6 & 1\\
\end{array}
$$
\qquad
$$
\begin{array}{cccccccc}
\end{array}
$$
\caption{\label{b=-1-2}The Motzkin triangle.}
\end{figure}

The bijection between trees and Motzkin paths is as follows: let $P$
be a Motzkin path. Let $0\!=\!i_1\!<\!\cdots\!<\!i_k\!=\!n$ be the sequence of
abscissas of integer points $(i,0)$ belonging to $P$ (also called the returns
to zero of $P$). Denote by $P_j$ the part of $P$ between $(i_j,0)$ and
$(i_{j+1},0)$.
Note that those elements have no non-trivial returns to zero.

Then the tree corresponding to $P$ is built in the following recursive way:
put $k$ at the root of the tree (meaning $d_k$).
If $P_j=(1,0)$ then put $c$ as the $j$-th son of the root. Else, $P_j$ is of
the form $P_j=(1,1)Q_j(1,-1)$. Then insert $Q_j$ recursively as the $j$-th son
of the root.

\begin{figure}[ht]
\begin{center}
\leavevmode
\epsfxsize =14cm
\epsffile{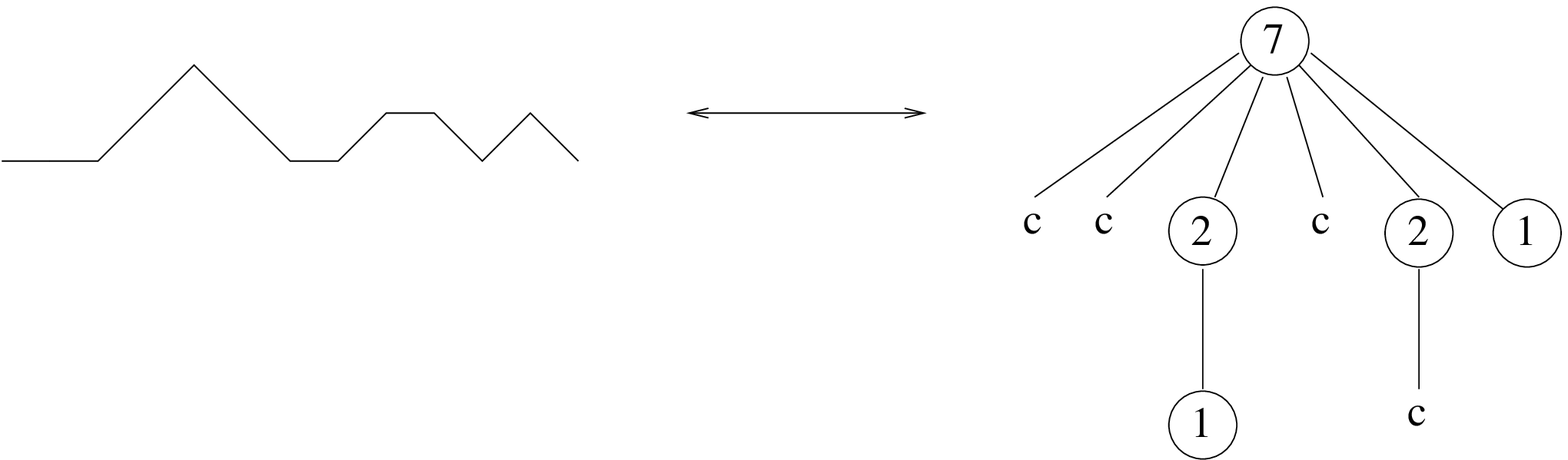}
\end{center}
\caption{\label{ex-bij-Motz}Example of the bijection between Motzkin paths and
trees.}
\end{figure}

Figure~\ref{ex-bij-Motz} presents an example of the bijection.

\subsection{Combinatorial sequences}

The generating functions of row-sums of the triangles are obtained by setting
$d_i=t^i$ in Equation (\ref{general-g}):
\begin{equation}
g = 1 + \frac{tg^{b+1}}{1-tg}.
\end{equation}
For $b=0$, one recognizes the quadratic equation satisfied by generating
series of Catalan numbers and for $b=1$ the quadratic equation of small
Schr\"oder numbers.
For $b=-1$, $g(-t)$ satisfies the quadratic equation for the Motzkin numbers.



\footnotesize
\begin{thebibliography}{aa}
%
\bibitem{BFK}{\sc C. Brouder, A. Frabetti} and {\sc C. Krattenthaler},
Non-commutative Hopf algebra of formal diffeomorphisms,
QA/0406117, 2004 - arxiv.org 
%
%
\bibitem{NCSF6}{\sc  G. Duchamp, F. Hivert}, and {\sc J.-Y. Thibon},
{\it Noncommutative symmetric functions VI: free quasi-symmetric functions and
related algebras},
Internat. J. Alg. Comput. {\bf 12} (2002), 671--717.
%
\bibitem{NCSF1}{\sc I.M. Gelfand, D. Krob, A. Lascoux, B. Leclerc,
V.~S. Retakh}, and {\sc J.-Y. Thibon},
{\it Noncommutative symmetric functions},
Adv. in Math. {\bf 112} (1995), 218--348.
%
%
\bibitem{Ges}{\sc I. Gessel},
{\it Noncommutative Generalization and $q$-analog of the Lagrange Inversion
Formula},
Trans. Amer. Math. Soc. {\bf 257} (1980), no. 2, 455--482.
%
\bibitem{GKP}
{\sc T. L. Graham, D. E. Knuth} and {\sc O. Patashnik},
{\it Concrete Mathematics}. Addison-Wesley, Reading, Mass., 1989;
2nd Ed. 1994.
%
\bibitem{Hai1}{\sc M. Haiman}, {\it 
Conjectures on the quotient ring by diagonal invariants},
J. Algebraic Combin. {\bf 3} (1994), 17--36.
%
\bibitem{Hiv} {\sc F. Hivert},
{\it Combinatoire des fonctions quasi-sym\'etriques},
Th\`ese de Doctorat, Marne-La-Vall\'ee, 1999.
%
\bibitem{HNT}{\sc F. Hivert, J.-C. Novelli}, {\sc J.-Y. Thibon},
{\it Commutative Hopf algebras of permutations and trees},
preprint math.CO/0502456.
%
\bibitem{HT}{\sc F. Hivert} and {\sc N. M. Thi\'ery},
{\it Representation theories of some towers of algebras related
to the symmetric groups and their Hecke algebras},
preprint (2005), submitted to FPSAC'06.
%
\bibitem{HTm} {\sc F. Hivert} and {\sc N. Thi\'ery},
{\it MuPAD-Combinat, an open-source package for research in algebraic
combinatorics},
S\'em. Lothar.  Combin. {\bf 51}  (2004), 70p. (electronic).
%
\bibitem{NCSF2}{\sc D. Krob, B. Leclerc} and {\sc J.-Y. Thibon},
{\it Noncommutative symmetric functions II: Transformations of alphabets},
Intern. J. Alg. Comput. {\bf 7} (1997), 181--264.
%
\bibitem{NCSF4}{\sc D. Krob} and {\sc J.-Y. Thibon},
{\it Noncommutative symmetric functions IV: Quantum linear groups and Hecke
algebras at $q=0$},
J. Alg. Comb. {\bf 6} (1997), 339--376.
%
\bibitem{KW}{\sc A. G. Konheim} and {\sc B. Weiss}, {\it
An occupancy discipline and applications},
SIAM J. Appl. Math. {\bf 14} (1966), 1266--1274.
%
%
\bibitem{KY}{\sc J. P. S. Kung} and {\sc C. Yan}, {\it
Gon\u carov polynomials and parking functions},
J. Combin. Theory A {\bf 102} (2003), 16--37. 
%
\bibitem{Las}{\sc A. Lascoux}, {\it
Symmetric functions and combinatorial operators on polynomials},
 CBMS Regional
Conference Series in Mathematics {\bf 99}, 
American Math. Soc.,
Providence, RI, 2003; xii+268 pp.
%
\bibitem{Len}{\sc C.  Lenart},
{\it Lagrange inversion and Schur functions},
J. Algebraic Combin. {\bf 11} (2000), 1, 69--78.
%
\bibitem{Mcd}{\sc I.G. Macdonald},
{\it Symmetric functions and Hall polynomials},
2nd ed., Oxford University Press, 1995.
%
%
\bibitem{NT1}{\sc J.-C Novelli} and {\sc J.-Y. Thibon},
{\it A Hopf algebra of parking functions},
Proc. FPSAC/SFCA 2004, Vancouver (electronic).
%
%
\bibitem{PPR}{\sc I. Pak, A. Postnikov}, and {\sc V.~S.  Retakh},
{\it Noncommutative Lagrange Theorem and Inversion Polynomials},
preprint, 1995, available at
{\tt http://www-math.mit.edu/${}^\sim$pak/research.html}.
%
\bibitem{Slo} {\sc N.J.A. Sloane},
{\it The On-Line Encyclopedia of Integer Sequences},\\
\verb+http://www.research.att.com/~njas/sequences/+
%
\bibitem{PS} {\sc R. P. Stanley} and {\sc J. Pitman}, {\it
A Polytope Related to Empirical Distributions, Plane Trees, Parking Functions,
and the Associahedron},
Discrete Comput. Geom.  {\bf 27}  (2002),  603--634.
%
\bibitem{Stan2}{\sc R. P. Stanley}, {\it Enumerative combinatorics},
vol. 2, Cambridge University Press, 1999.
%
\bibitem{Yan}{\sc C. H. Yan}, {\it Generalized parking functions, tree
inversions and multicolored graphs},
Adv. Appl. Math. {\bf 27} (2001), 641--670.
%
\end{thebibliography}
\end{document}